\documentclass[11pt]{article}%
\usepackage{amsmath}
\usepackage{amsfonts}
\usepackage{amssymb}
\usepackage{graphicx}
\usepackage{theorem}
\usepackage{makeidx}
\usepackage{verbatim}
\usepackage{color}%
\providecommand{\U}[1]{\protect\rule{.1in}{.1in}}
\newtheorem{thm}{Theorem}[section]

\newtheorem{pr}[thm]{Proposition}
\newtheorem{df}[thm]{Definition}
\newtheorem{rmk}[thm]{Remark}

{\theorembodyfont{\upshape}
\newtheorem{examp}[thm]{Example}
}
\numberwithin{equation}{section} 
\begin{document}

\title{Title: }

\begin{center}

\vspace*{1.3cm}

\textbf{DIRECTIONAL PARETO\ EFFICIENCY: CONCEPTS AND OPTIMALITY CONDITIONS}

\bigskip

by

\bigskip

TEODOR CHELMU\c{S}\footnote{{\small {Faculty of Mathematics, "Alexandru Ioan
Cuza" University,} \ {Bd. Carol I, nr. 11, 700506 -- Ia\c{s}i, Romania,
e-mail: \texttt{teo1chelmus@gmail.com}}}},
MARIUS\ DUREA\footnote{{\small {Faculty of Mathematics, "Alexandru Ioan Cuza"
University,} \ {Bd. Carol I, nr. 11, 700506 -- Ia\c{s}i, Romania, e-mail:
\texttt{durea@uaic.ro}}}, {\small and "Octav Mayer" Institute of Mathematics
of the Romanian Academy, Ia\c{s}i, Romania.}}, and ELENA-ANDREEA
FLOREA\footnote{{\small {Faculty of Mathematics, "Alexandru Ioan Cuza"
University,} \ {Bd. Carol I, nr. 11, 700506 -- Ia\c{s}i, Romania, e-mail:
\texttt{andreea\_acsinte@yahoo.com}}}}
\end{center}

\bigskip

\noindent{\small {\textbf{Abstract:}} We introduce and study a notion of
directional Pareto minimality with respect to a set that generalizes the
classical concept of Pareto efficiency. Then we give separate necessary and
sufficient conditions for the newly introduced efficiency and several
situations concerning the objective mapping and the constraints are
considered. In order to investigate different cases, we adapt some well-known
constructions of generalized differentiation and the connections with some
recent directional regularities come naturally into play. As a consequence,
several techniques from the study of genuine Pareto minima are considered in
our specific situation.}

\bigskip

\noindent{\small {\textbf{Keywords:}} {directional }Pareto minimality{ $\cdot$
optimality conditions $\cdot$ directional tangent cones} {$\cdot$\textbf{
}directional regularity}}

\bigskip

\noindent{\small {\textbf{Mathematics Subject Classification (2010): }}54C60
}$\cdot$ {\small 46G05 }$\cdot$ {\small 90C46}

\begin{center}

\end{center}

\section{Introduction and notation}

This paper has two main motivations. On one hand, we are aiming at continuing
the effort made by several authors in the last decade to investigate
directional phenomena in mathematical programming and, on the other hand, we
show the power of several tools related to directional regularities that have
been developed recently. For detailed accounts on these topics we refer the
reader to the following works and references therein: \cite{Mord&nam},
\cite{nam&zalinescu}, \cite{Popovici.tammer}, \cite{HGfrerer},
\cite{Thera&Huynh}.

In this work, inspired by some ideas coming in vector optimization problems
from location theory where some directions are privileged with respect to the
others, we present a notion of directional minimality for mappings and we
illustrate by examples its relevance even for the case of real-valued
functions. Then, we observe on the simplest case of real-valued functions of a
real variable that the natural necessary optimality conditions are given by
the Fermat Theorem at an endpoint of an interval. This gives us the impetus to
consider far-reaching generalization of this case, namely, problems where the
objective is a set-valued map and the constraint is given by means of an
inverse image of a cone through another set-valued map. For the study of this
general case, we introduce an adapted tangent cone, along with several
directional regularity properties of the involved maps, and this approach
allows us to derive necessary optimality conditions that, in turn, generalize
the prototype of Fermat Theorem at an endpoint of an interval. Furthermore, we
present as well optimality conditions in terms of tangent limiting cones and
coderivatives. Both on primal and dual spaces we have under consideration
several situations concerning the objective and constraint mappings with their
specific techniques of study, among which we mention generalized constraint
qualification conditions, Gerstewitz scalarization, openness vs. minimality
paradigm, Clarke penalization, extremal principle. Some results are dedicated
to the sufficient optimality conditions under convexity assumptions. Finally,
we consider as well the situation of minimality for sets and a brief
discussion of this concept reveals the similarities and the differences with
respect to the known situation of Pareto efficiency.

The paper is organized as follows. First of all, we introduce the notation we
use and then we present the concepts of directional minimality we study in
this work. The definitions of these notions along with some comparisons and
examples are the subjects of the second section. The main section of the paper
is the third one, and it deals with optimality conditions for the above
introduced concepts, being, in turn, divided into two subsections. Firstly, we
derive optimality conditions using tangent cones and to this aim we adapt a
classical concept of the Bouligand tangent cone and Bouligand derivative of a
set-valued map. Using some directional metric regularities, we get several
assertions concerning these objects and this allows us to present necessary
optimality conditions for a wide range of situations going from problems
governed by set-valued mappings having generalized inequalities constraints to
fully smooth constrained problems. Secondly, we deal with optimality
conditions using normal limiting cones and, again, we consider several types
of problems. In this process of getting necessary optimality conditions we
adapt several techniques from classical vector optimization. Moreover, some
generalized convex cases are considered in order to obtain sufficient
optimality conditions. The last section deals with Pareto directional minima
for sets. We emphasize the fact that even if the directional Pareto efficiency
appears naturally in the case of mappings, it can be considered as well for
sets and in this respect we present the corresponding concepts and we discuss
it by means of some examples and optimality conditions in terms of the
modified tangent cones. Several conclusions of this work are collected in a
short section that ends the paper. \bigskip

Throughout this paper, we assume that $X$, $Y$ and $Z$ are normed vector
spaces over the real field $\mathbb{R}$ and on a product of normed vector
spaces we consider the sum norm, unless otherwise stated. By $B\left(
x,\varepsilon\right)  $ we denote the open ball with center $x$ and radius
$\varepsilon>0$ and by $B_{X}$ the open unit ball of $X.$ In the same manner,
$D(x,\varepsilon)$ and $D_{X}$ denote the corresponding closed balls. The
symbol $S_{X}$ stands for the unit sphere of $X.$ By the symbol $X^{\ast}$ we
denote the topological dual of $X$, while $w^{\ast}$ stands for the
weak$^{\ast}$ topology on $X^{\ast}.$

Let $F:X\rightrightarrows Y$ be a set-valued map. As usual, the graph of $F$
is%
\[
\operatorname*{Gr}F:=\left\{  \left(  x,y\right)  \in X\times Y\mid y\in
F\left(  x\right)  \right\}  ,
\]
and the inverse of $F$ is the set-valued map $F^{-1}:Y\rightrightarrows X$
given by $\left(  y,x\right)  \in\operatorname*{Gr}F^{-1}$ iff $\left(
x,y\right)  \in\operatorname*{Gr}F$. Consider a nonempty subset $A$ of $X.$
Then the image of $A$ through $F$ is
\[
F\left(  A\right)  :=\left\{  y\in Y\mid\exists x\in A:\text{ }y\in F\left(
x\right)  \right\}
\]
and the distance function associated to $A$ is $d_{A}:X\rightarrow\mathbb{R}$
given by%
\[
d_{A}\left(  x\right)  =d(x,A):=\inf\limits_{a\in A}\left\Vert x-a\right\Vert
.
\]
The topological interior, topological closure, the convex hull and conic hull
of $A$ are denoted, respectively, by $\operatorname*{int}A,$
$\operatorname*{cl}A,$ $\operatorname*{conv}A,$ $\operatorname*{cone}A$. The
negative polar of $A$ is
\[
A^{-}:=\left\{  x^{\ast}\in X^{\ast}\mid x^{\ast}\left(  a\right)
\leq0,\forall a\in A\right\}  .
\]

\section{The concepts under study}

Let $K\subset Y$ be a proper (that is, $K\neq\{0\},$ $K\neq Y$) convex cone
(we do not suppose that $K$ is pointed, in general). For such a cone, its
positive dual cone is
\[
K^{+}:=\left\{  y^{\ast}\in Y^{\ast}\mid y^{\ast}\left(  y\right)
\geq0,\forall y\in K\right\}  .
\]

Take $F:X\rightrightarrows Y$ as a set-valued mapping,$\ $and let us consider
the following geometrically constrained optimization problem with
multifunctions:
\[
(P)\hspace{0.4cm}\operatorname*{minimize}\text{ }F(x),\text{ subject to }x\in
A,
\]
where $A\subset X$ is a closed nonempty set.

Usually, the minimality is understood in the Pareto sense given by the next definition.

\begin{df}
A point $(\overline{x},\overline{y})\in\operatorname*{Gr}F\cap(A\times Y)$ is
a local Pareto minimum point for $F$ on $A$ if there exists a neighborhood $U$
of $\overline{x}$ such that
\begin{equation}
\left(  F(U\cap A)-\overline{y}\right)  \cap-K\subset K. \label{P_rest}%
\end{equation}

\end{df}

The vectorial notion described by (\ref{P_rest}) covers as well the situation
where $f$ is a function (in which case $\overline{y}=f(\overline{x})$ will not
be mentioned) and the situation of classical local minima in scalar case (in
which case we drop the label "Pareto"). If $K$ is pointed (that is,
$K\cap-K=\{0\}$) then (\ref{P_rest}) reduces to
\[
\left(  F(U\cap A)-\overline{y}\right)  \cap-K\subset\{0\}.
\]

\begin{df}
If $\operatorname*{int}K\neq\emptyset,$ the point $(\overline{x},\overline
{y})\in\operatorname*{Gr}F\cap(A\times Y)$ is a local weak Pareto minimum
point for $F$ on $A$ if there exists a neighborhood $U$ of $\overline{x}$ such
that
\[
\left(  F(U\cap A)-\overline{y}\right)  \cap-\operatorname*{int}K=\emptyset.
\]

\end{df}

Let $L\subset S_{X}$ be a nonempty closed set. Then it is not difficult to see
that $\operatorname*{cone}L$ is closed as well. Indeed, let us consider a
sequence $(u_{n})\subset\operatorname*{cone}L$ converging towards $u\in X.$ We
have to show that $u\in\operatorname*{cone}L.$ The case $u=0\in
\operatorname*{cone}L$ is clear$.$ Otherwise, there are some sequences
$(t_{n})\subset(0,\infty)$ and $(\ell_{n})\subset L$ such that $u_{n}%
=t_{n}\ell_{n}$ for every $n.$ If $(t_{n})\rightarrow0$ (on a subsequence),
the boundedness of $(\ell_{n})$ leads to $u=0,$ a situation avoided at this
stage. If $(t_{n})$ is unbounded, then again the relation $\left\Vert
t_{n}\ell_{n}\right\Vert \rightarrow\left\Vert u\right\Vert $ leads to a
contradiction. So, on a subsequence, $(t_{n})\rightarrow t>0$ which means, by
the closedness of $L,$ that $\ell_{n}=t_{n}^{-1}t_{n}\ell_{n}\rightarrow
t^{-1}u\in L$, therefore $u\in\operatorname*{cone}L,$ as claimed.

\bigskip

The main purpose of this paper is to introduce and to study the following concept.

\begin{df}
\label{df_dir_min_mf}One says that $(\overline{x},\overline{y})\in
\operatorname*{Gr}F\cap(A\times Y)$ is a local directional Pareto minimum
point for $F$ on $A$ with respect to (the set of directions) $L$ if there
exists a neighborhood $U$ of $\overline{x}$ such that
\begin{equation}
\left(  F(U\cap A\cap\left(  \overline{x}+\operatorname*{cone}L\right)
)-\overline{y}\right)  \cap-K\subset K. \label{P_dir}%
\end{equation}

\end{df}

If one compares this relation to (\ref{P_rest}), then one observes that this
concept corresponds to the situation where the restriction has the special
form (depending on the reference point) $A\cap\left(  \overline{x}%
+\operatorname*{cone}L\right)  .$ Of course, when $A=X$ in (\ref{P_dir}) then
one says that $(\overline{x},\overline{y})\in\operatorname*{Gr}F$ is a local
directional Pareto minimum point for $F$ with respect to $L.$ Now, the concept
of local directional Pareto maximum is obtained in an obvious way.

If $\operatorname*{int}K\neq\emptyset,$ one defines as well the weak
counterpart of the above notion.

\begin{df}
One says that $(\overline{x},\overline{y})\in\operatorname*{Gr}F\cap(A\times
Y)$ is a local weak directional Pareto minimum point for $F$ on $A$ with
respect to (the set of directions) $L$ if there exists a neighborhood $U$ of
$\overline{x}$ such that
\[
\left(  F(U\cap A\cap\left(  \overline{x}+\operatorname*{cone}L\right)
)-\overline{y}\right)  \cap-\operatorname*{int}K=\emptyset.
\]

\end{df}

In all these notions, if one takes $U=X,$ then we get the corresponding global concepts.

\begin{rmk}
If $L_{1},L_{2}\subset S_{X}$ are nonempty closed subsets such that
$L_{1}\subset L_{2},$ then a local directional Pareto minimum point for $F$
with respect to $L_{2}$ is a local directional Pareto minimum point for $F$
with respect to $L_{1}.$
\end{rmk}

It is obvious that (\ref{P_rest}) implies (\ref{P_dir}), but the converse is
not true. To justify the latter affirmation, let us consider the following
simple scalar example (when the output space is $\mathbb{R}$ we always
consider $K:=[0,\infty)$)$.$

\begin{examp}
Let $f:\mathbb{R\rightarrow R}$ be a strictly increasing function. Then every
$\overline{x}\in\mathbb{R}$ is local directional minimum for $f$ with respect
to $L:=\{+1\},$ but it is not a local minimum for $f.$
\end{examp}

Moreover, the minimality concept introduced here covers some interesting
situations described by the next examples.

\begin{examp}
Let $f:\mathbb{R}^{2}\mathbb{\rightarrow R}$ be given by $f(x,y)=x^{2}-y^{2}.$
It is well known that $(0,0)$ is a critical saddle point, whence it is not a
minimum point. However, it is a directional minimum point for $f$ with respect
to $L=\{-1,1\}\times\{0\}$ since for every $(x,y)\in(0,0)+\operatorname*{cone}%
L=\mathbb{R\times\{}0\},$ one has $f(x,y)\geq f(0,0).$ Similarly, $(0,0)$ is a
directional maximum point for $f$ with respect to $L=\{0\}\times\{-1,1\}.$
\end{examp}

\begin{examp}
Let $f:\mathbb{R}^{2}\mathbb{\rightarrow R}$ be given by $f(x,y)=x^{2}-y^{3}.$
Again, $(0,0)$ is a critical saddle point. It is now easy to see that it is,
however, a directional minimum point for $f$ with respect to $L=\{-1,1\}\times
\{0\}$ and to respect to $L=\{0\}\times\{-1\}.$
\end{examp}

The next example emphasizes that there are points which are not directional
minima with respect to any nonempty closed set $L\subset S_{X}.$ This applies
also for critical points of smooth functions.

\begin{examp}
Let $f:\mathbb{R\rightarrow R}$ be given by%
\[
f(x)=\left\{
\begin{array}
[c]{l}%
\sin\frac{1}{x},\text{ if }x\neq0\\
0,\text{ if }x=0
\end{array}
\right.  .
\]
Then $\overline{x}=0$ is not directional minimum for $f$ neither for
$L:=\{-1\},$ nor for $L:=\{+1\}.$ In the same manner, $f:\mathbb{R\rightarrow
R}$ given by%
\[
f(x)=\left\{
\begin{array}
[c]{l}%
x^{3}\sin\frac{1}{x},\text{ if }x\neq0\\
0,\text{ if }x=0
\end{array}
\right.
\]
is differentiable at $\overline{x}=0,$ $f^{\prime}(\overline{x})=0,$ but
$\overline{x}$ is not a directional minimum for $f.$
\end{examp}

The next example underlines the idea that for every prescribed set of
directions one can define functions that achieve directional minimum with
respect to the given set.

\begin{examp}
Let $0<\theta_{1}<\theta_{2}<\frac{\pi}{2}$ and $L:=\{(\cos\theta,\sin
\theta)\mid\theta_{1}\leq\theta\leq\theta_{2}\}.$ Consider $f:\mathbb{R}%
^{2}\mathbb{\rightarrow R}$ be given by%
\[
f(x,y)=\left\{
\begin{array}
[c]{l}%
\left(  \theta_{2}-\arctan\frac{y}{x}\right)  \left(  \arctan\frac{y}%
{x}-\theta_{1}\right)  ,\text{ if }x\neq0\text{ and }(x>0\text{ or }y\geq0)\\
0,\text{ if }x=0\\
-1,\text{ if }x<0\text{ and }y<0
\end{array}
\right.  .
\]
Then it is not difficult to see that $(0,0)$ is directional minimum for $f$
with respect to $L$.
\end{examp}

Using these basic examples of scalar-valued functions, we are able to easily
build examples for vector-valued maps. Here are two such examples.

\begin{examp}
Consider $f:\mathbb{R}^{2}\mathbb{\rightarrow R}^{2}$ be given by
$f(x,y)=(x^{2}-y^{2},x^{2}-y^{3}).$ Consider $K:=\mathbb{R}_{+}^{2}.$ Then
$(0,0)$ is a directional minimum for $f$ with respect to $L:=\{(1,0)\}\subset
S_{\mathbb{R}^{2}}.$
\end{examp}

\begin{examp}
Let $f:\mathbb{R\rightarrow R}^{2}$ be given by $f\left(  x\right)  =\left(
2x,x\right)  $ and $K=\operatorname*{cone}\operatorname*{conv}\left\{  \left(
1,0\right)  ,\left(  1,1\right)  \right\}  .$ It is easy to see that
$\overline{x}:=0$ is a directional minimum for $f$ with respect to
$L:=\left\{  +1\right\}  ,$ but $\overline{x}$ is not a local Pareto minimum
point for $f.$
\end{examp}

The concepts introduced in this section are studied in the sequel from the
point of view of optimality conditions.

\section{Optimality conditions for directional minima}

In order to start with the necessary optimality conditions for directional
minima, let us to observe that the obvious prototype for such an investigation
is the Fermat Theorem for derivable real-valued functions with one variable at
interval endpoints: if $f:[a,b]\rightarrow\mathbb{R}$ is a function for which
$a$ is local minimum point (that is, a directional minimum with respect to
$L:=\{+1\}$), and $f$ is derivable at $a,$ then $f^{\prime}(a)\geq0,$ and,
similarly, if $b$ is a minimum point for $f$ (that is, a directional minimum
with respect to $L:=\{-1\}$)$,$ and $f$ is derivable at $b,$ then $f^{\prime
}(b)\leq0.$

We approach this issue from two points of view, namely, making use of tangent
cones (which are objects of generalized differentiation on primal spaces) and
of normal cones (constructions that are defined on dual spaces).

\subsection{Optimality conditions using tangent cones}

Let us consider now several concepts that will help us in studying optimality
conditions for the directional minima.

\begin{df}
Let $A\subset X$ be a nonempty set and $L\subset S_{X}$ be a nonempty closed
set. Then the Bouligand tangent cone to $A$ at $\overline{x}\in A$ with
respect to $L$ is the set
\[
T_{B}^{L}(A,\overline{x}):=\left\{  u\in X\mid\exists(u_{n})\overset
{\operatorname*{cone}L}{\longrightarrow}u,\exists(t_{n})\overset{(0,\infty
)}{\longrightarrow}0\text{ such that for all }n,\text{ }\overline{x}%
+t_{n}u_{n}\in A\right\}  ,
\]
where $(u_{n})\overset{\operatorname*{cone}L}{\longrightarrow}u$ means
$(u_{n})\longrightarrow u$ and $(u_{n})\subset\operatorname*{cone}L$, and
similarly for $(t_{n})\overset{(0,\infty)}{\longrightarrow}0.$
\end{df}

Obviously, this is a adaptation of the concept of Bouligand tangent cone to
$A$ at $\overline{x}$ defined as%
\[
T_{B}(A,\overline{x}):=\left\{  u\in X\mid\exists(u_{n})\rightarrow
u,\exists(t_{n})\overset{(0,\infty)}{\longrightarrow}0\text{ such that for all
}n,\text{ }\overline{x}+t_{n}u_{n}\in A\right\}  .
\]

Some remarks are in order.

\begin{rmk}
As the usual Bouligand tangent cone, the set $T_{B}^{L}(A,\overline{x})$ is a
closed cone: the proof of this assertion can be made directly as for the
classical concept (see \cite{AF}) or by observing that%
\[
T_{B}^{L}(A,\overline{x})=T_{B}(A\cap\left(  \overline{x}+\operatorname*{cone}%
L\right)  ,\overline{x}).
\]
In view of the fact that $\operatorname*{cone}L$ is closed, one has that
$T_{B}^{L}(A,\overline{x})\subset\operatorname*{cone}L.$ Moreover,
\[
T_{B}^{L}(A,\overline{x})\subset T_{B}(A,\overline{x})\cap T_{B}(\overline
{x}+\operatorname*{cone}L,\overline{x})=T_{B}(A,\overline{x})\cap
\operatorname*{cone}L.
\]
However, the inclusion above does not hold as equality, in general. To see
this, consider the set $A\subset X:=\mathbb{R}^{2}$ as the plane domain
bounded by the curve (the cardioid), which has the parametric representation%
\[
\left\{
\begin{array}
[c]{l}%
x=-2\cos t+\cos2t+1\\
y=2\sin t-\sin2t
\end{array}
,\text{ }t\in\lbrack0,2\pi]\right.  ,
\]
$\overline{x}:=(0,0),$ $L:=\{\left(  -1,0\right)  \}$ and observe that
$T_{B}(A,\overline{x})=X$ and $T_{B}(A\cap\left(  \overline{x}%
+\operatorname*{cone}L\right)  ,\overline{x})=\{\overline{x}\}.$

Another useful and easy-to-see inclusion is%
\[
\operatorname*{cl}\left(  T_{B}(A,\overline{x})\cap\operatorname*{int}%
\operatorname*{cone}L\right)  \subset T_{B}^{L}(A,\overline{x}).
\]

\end{rmk}

\begin{df}
Let $F:X\rightrightarrows Y$ be a set-valued map, $(\overline{x},\overline
{y})\in\operatorname*{Gr}F$ and $L\subset S_{X},$ $M\subset S_{Y}$ be nonempty
closed sets. The Bouligand derivative of $F$ at $(\overline{x},\overline{y})$
with respect to $L$ and $M$ is the set-valued map $D_{B}^{L,M}F(\overline
{x},\overline{y}):X\rightrightarrows Y$ defined by the relation $v\in
D_{B}^{L,M}F(\overline{x},\overline{y})(u)$ iff there are $(u_{n}%
)\overset{\operatorname*{cone}L}{\longrightarrow}u,(v_{n})\overset
{\operatorname*{cone}M}{\longrightarrow}v,(t_{n})\overset{(0,\infty
)}{\longrightarrow}0$ such that for all $n,$
\[
\overline{y}+t_{n}v_{n}\in F(\overline{x}+t_{n}v_{n}).
\]

\end{df}

Clearly,%
\[
\operatorname*{Gr}D_{B}^{L,M}F(\overline{x},\overline{y})\subset
\operatorname*{cone}L\times\operatorname*{cone}M.
\]

Again, this is an adaptation of the well-known Bouligand derivative of $F$ at
$(\overline{x},\overline{y}),$ which is the set-valued map $D_{B}%
F(\overline{x},\overline{y}):X\rightrightarrows Y$ defined by%
\[
\operatorname*{Gr}D_{B}F(\overline{x},\overline{y}):=T_{B}\left(
\operatorname*{Gr}F,(\overline{x},\overline{y})\right)  .
\]

Other derivability objects in primal spaces that can be adapted in directional
setting in a similar manner are the Ursescu (adjacent) tangent cone and the
Ursescu (adjacent) derivative (see \cite{D2004}), and the Dini lower
derivative of $F$ at $(\overline{x},\overline{y}),$ which is the multifunction
$D_{D}F(\overline{x},\overline{y})$ from $X$ into $Y$ given, for every $u\in
X,$ by
\begin{align*}
D_{D}F(\overline{x},\overline{y})(u)=\{  &  v\in Y\mid\forall(t_{n}%
)\overset{(0,\infty)}{\longrightarrow}0,\forall(u_{n})\rightarrow
u,\exists(v_{n})\rightarrow v,\\
\forall n  &  \in\mathbb{N},\overline{{y}}+t_{n}v_{n}\in F(\overline{x}%
+t_{n}u_{n})\}.
\end{align*}
When $F:=f$ is a single-valued map, for simplicity, we write $D_{B}%
^{L,M}f(\overline{x})$ for $D_{B}^{L,M}f(\overline{x},\overline{y})$, and
similarly for $D_{D}$.

We present now the first result of this work.

\begin{pr}
\label{prop_weak}In the above notation, if $\operatorname*{int}K\neq\emptyset$
and $(\overline{x},\overline{y})\in\operatorname*{Gr}F$ is a local weak
directional Pareto minimum point for $F$ on $A$ with respect to $L$ then%
\[
D_{D}F(\overline{x},\overline{y})(T_{B}^{L}(A,\overline{x}))\cap
-\operatorname*{int}K=\emptyset.
\]
Moreover, if $A=X,$ then
\[
D_{B}^{L,S_{Y}}F(\overline{x},\overline{y})(X)\cap-\operatorname*{int}%
K=\emptyset.
\]

\end{pr}

\noindent\textbf{Proof. }We prove only the second part, since the first part,
on one hand, is similar, and, on the other hand, it follows from the
definitions and \cite[Theorem 3.1]{D2004}. Take $u\in X.$ If $u\notin
\operatorname*{cone}L$ then $D_{B}^{L,S_{Y}}F(\overline{x},\overline
{y})(u)=\emptyset$ and there is nothing to prove. If $u\in\operatorname*{cone}%
L,$ suppose, by way of contradiction, that there is $k\in-\operatorname*{int}%
K$ such that%
\[
k\in D_{B}^{L,S_{Y}}F(\overline{x},\overline{y})(u).
\]
According to the definition of $D_{B}^{L,S_{Y}}F(\overline{x},\overline{y}),$
this means that there exist $(t_{n})\overset{(0,\infty)}{\longrightarrow}0,$
$(u_{n})\overset{\operatorname*{cone}L}{\longrightarrow}u,$ $(k_{n}%
)\rightarrow k$ such that for all $n,$%
\[
\overline{y}+t_{n}k_{n}\in F(\overline{x}+t_{n}u_{n}),
\]
that is,
\[
t_{n}k_{n}\in F(\overline{x}+t_{n}u_{n})-\overline{y}.
\]
But, for $n$ large enough, $\overline{x}+t_{n}u_{n}$ is close enough to
$\overline{x}$ and belongs as well to $\overline{x}+\operatorname*{cone}L.$
Then, for such $n,$ taking into account the minimality of $(\overline
{x},\overline{y}),$ one gets $t_{n}k_{n}\notin-\operatorname*{int}K$ which
contradicts the fact that $k_{n}\rightarrow k\in-\operatorname*{int}K.$%
\hfill$\square$

\bigskip

In \cite{DPS}, by means of a special type of minimal time function, several
directional regularity properties for set-valued maps are introduced and
studied. In order to further investigate the directional minima we need to
briefly point out the main aspects concerning the minimal time function and
some related directional metric regularity.

Consider $\emptyset\neq L\subset S_{X}$ and $\emptyset\neq\Omega\subset X.$
Then the function%
\begin{align}
T_{L}(x,\Omega)  &  :=\inf\left\{  t\geq0\mid\exists u\in L:x+tu\in
\Omega\right\} \label{min_t_M}\\
&  =\inf\left\{  t\geq0\mid(x+tL)\cap\Omega\neq\emptyset\right\} \nonumber
\end{align}
is called the directional minimal time function with respect to $L.$

Remark that, if $L=S_{X},$ then $T_{L}(\cdot,\Omega)=d(\cdot,\Omega)$.
Moreover, we add the convention that $T_{L}(x,\emptyset)=\infty$ for every $x$
and we denote in what follows $T_{L}(x,\left\{  u\right\}  )$ by $T_{L}(x,u).$
Obviously, $T_{L}(x,u)<+\infty$ is equivalent to $T_{L}(x,u)=\left\Vert
u-x\right\Vert $ and $u-x\in\operatorname*{cone}L.$

Let $F:X\rightrightarrows Y$ be a set-valued mapping and $(\overline
{x},\overline{y})\in\operatorname*{Gr}F,$ $\emptyset\neq L\subset S_{X},$
$\emptyset\neq M\subset S_{Y}.$

What we need in the sequel is the following concept of directional calmness.
One says that $F$ is directionally calm at $(\overline{x},\overline{y})$ with
respect to $L$ and $M$ if there are $\alpha>0$ and some neighborhoods $U\ $of
$\overline{x}$ and $V$ of $\overline{y}$ such that for every $x\in U,$%

\begin{equation}
\sup_{y\in F(x)\cap V}T_{M}(y,F(\overline{x}))\leq\alpha T_{L}(\overline
{x},x). \label{calm}%
\end{equation}

We use the convention $\sup\limits_{x\in\emptyset}T_{L}\left(  x,\Omega
\right)  :=0$ for every nonempty set $\Omega\subset X$.

As usual (see \cite[Section 3H]{DontRock2009b}), for a calmness concept for
$F$, it is natural to have a metric subregularity notion such that the former
property for $F^{-1}$ to be equivalent to the latter property for $F.$ In our
setting, this corresponding concept reads as follows: one says that $F$ is
directionally metric subregular at $(\overline{x},\overline{y})$ with respect
to $L$ and $M$ if there exist $\alpha>0$ and some neighborhoods $U\ $of
$\overline{x}$ and $V$ of $\overline{y}$ such that for every $x\in U,$%
\begin{equation}
T_{L}(x,F^{-1}(\overline{y}))\leq\alpha T_{M}(\overline{y},F(x)\cap V).
\label{metric_subreg}%
\end{equation}
The expected equivalence is described in the following result.

\begin{pr}
\label{prop_echiv_mr_c}The set-valued map $F$ is directionally metric
subregular at $(\overline{x},\overline{y})$ with respect to $L$ and $M$ iff
$F^{-1}$ is directionally calm at $(\overline{y},\overline{x})$ with respect
to $M$ and $L.$
\end{pr}

\noindent\textbf{Proof.} Suppose first that $F$ is directionally metric
subregular at $(\overline{x},\overline{y})$ with respect to $L$ and $M.$ Then,
there exist $\alpha>0,$ $U\in\mathcal{V}\left(  \overline{x}\right)  $ and
$V\in\mathcal{V}\left(  \overline{y}\right)  $ such that for every $x\in U$
relation (\ref{metric_subreg}) holds. Let $y\in V.$ If $T_{M}(\overline
{y},y)=+\infty,$ there is nothing to prove. Suppose that $T_{M}(\overline
{y},y)<+\infty,$ which means that $y-\overline{y}\in\operatorname*{cone}M.$
Consider $x\in U$ with $y\in F(x),$ i.e., $x\in F^{-1}(y)\cap U.$ Then, by
hypothesis,
\[
T_{L}(x,F^{-1}(\overline{y}))\leq\alpha T_{M}(\overline{y},F(x)\cap
V)\leq\alpha T_{M}(\overline{y},y),
\]
so,%
\[
\sup_{x\in F^{-1}(y)\cap U}T_{L}(x,F^{-1}(\overline{y}))\leq\alpha
T_{M}(\overline{y},y),
\]
for all $y\in V,$ whence the conclusion.

For the converse, suppose that $F^{-1}$ is directionally calm at
$(\overline{y},\overline{x})$ with respect to $M$ and $L.$ Therefore, there
exist $\alpha>0,$ $U\in\mathcal{V}\left(  \overline{x}\right)  $ and
$V\in\mathcal{V}\left(  \overline{y}\right)  $ such that for every $y\in V$%
\[
\sup_{x\in F^{-1}(y)\cap U}T_{L}(x,F^{-1}(\overline{y}))\leq\alpha
T_{M}(\overline{y},y).
\]
Take $x\in U$. Again, if $T_{M}(\overline{y},F(x)\cap V)=+\infty,$ the desired
inequality holds. Suppose that $T_{M}(\overline{y},F(x)\cap V)<+\infty,$ which
means that for any $\varepsilon>0$ there exist $u_{\varepsilon}\in M$ and
$y_{\varepsilon}\in F(x)\cap V$ such that%
\[
\overline{y}+\left(  T_{M}(\overline{y},F(x)\cap V)+\varepsilon\right)
u_{\varepsilon}=y_{\varepsilon}.
\]
Therefore, $y_{\varepsilon}-\overline{y}\in\operatorname*{cone}M,$ $x\in
F^{-1}(y_{\varepsilon})\cap U$ and from the hypothesis,
\[
T_{L}(x,F^{-1}(\overline{y}))\leq\alpha T_{M}(\overline{y},y_{\varepsilon
})=\alpha\left\Vert y_{\varepsilon}-\overline{y}\right\Vert \leq\alpha\left(
T_{M}(\overline{y},F(x)\cap V)+\varepsilon\right)  .
\]
Passing to the limit as $\varepsilon\rightarrow0$ we get the conclusion.\hfill
$\square$

\bigskip

Now, we use the directional calmness for getting an evaluation of the
directional Bouligand tangent cone to a value of a set-valued mapping in terms
of the image of $0$ through the directional Bouligand derivative of the same application.

\begin{pr}
\label{prop_calm}Let $F:X\rightrightarrows Y$ be a set-valued mapping,
$(\overline{x},\overline{y})\in\operatorname*{Gr}F$, and $\emptyset\neq
L\subset S_{X},$ $\emptyset\neq M\subset S_{Y}$ be closed sets. Then
\[
T_{B}^{M}(F\left(  \overline{x}\right)  ,\overline{y})\subset D_{B}%
^{L,M}F(\overline{x},\overline{y})(0).
\]
Moreover, if $F\ $is directionally calm at $(\overline{x},\overline{y})$ with
respect to $L$ and $M,$ and $\operatorname*{cone}M$ is convex, then the
equality holds.
\end{pr}

\noindent\textbf{Proof.} Take $v\in T_{B}^{M}(F\left(  \overline{x}\right)
,\overline{y}).$ According to the definition, there are $(v_{n})\overset
{\operatorname*{cone}M}{\longrightarrow}v,(t_{n})\overset{(0,\infty
)}{\longrightarrow}0$ such that for all $n,$
\[
\overline{y}+t_{n}v_{n}\in F(\overline{x})=F(\overline{x}+t_{n}\cdot0),
\]
which clearly implies that $v\in D_{B}^{L,M}F(\overline{x},\overline{y})(0).$

For the opposite inclusion, take $v\in D_{B}^{L,M}F(\overline{x},\overline
{y})(0)$ meaning that there are $(u_{n})\overset{\operatorname*{cone}%
L}{\longrightarrow}0,(v_{n})\overset{\operatorname*{cone}M}{\longrightarrow
}v,(t_{n})\overset{(0,\infty)}{\longrightarrow}0$ such that for all $n,$
\[
\overline{y}+t_{n}v_{n}\in F(\overline{x}+t_{n}u_{n}).
\]
But, the assumed calmness of $F$ and the fact that $\left(  t_{n}u_{n}\right)
\subset\operatorname*{cone}L,$ mean that, for a positive $\alpha$ and for all
$n$ large enough,
\[
T_{M}(\overline{y}+t_{n}v_{n},F(\overline{x}))\leq\alpha T_{L}(\overline
{x},\overline{x}+t_{n}u_{n})=\alpha t_{n}\left\Vert u_{n}\right\Vert ,
\]
that is
\[
\inf\{\tau\geq0\mid\exists w_{n}\in M\text{ such that for all }n,\text{
}\overline{y}+t_{n}v_{n}+\tau w_{n}\in F(\overline{x})\}\leq\alpha
t_{n}\left\Vert u_{n}\right\Vert .
\]
Therefore, for every $n$ (large enough) there are $w_{n}\in M$ and $\tau
_{n}\geq0$ such that $\beta_{n}:=\overline{y}+t_{n}v_{n}+\tau_{n}w_{n}\in
F(\overline{x})$ and $\tau_{n}<\alpha t_{n}\left\Vert u_{n}\right\Vert
+t_{n}^{2}.$ So, for every $n,$
\[
\left\Vert \beta_{n}-\left(  \overline{y}+t_{n}v_{n}\right)  \right\Vert
=\tau_{n}<\alpha t_{n}\left\Vert u_{n}\right\Vert +t_{n}^{2},
\]
whence%
\[
\left\Vert \frac{1}{t_{n}}(\beta_{n}-\overline{y})-v_{n}\right\Vert
<\alpha\left\Vert u_{n}\right\Vert +t_{n},
\]
which gives%
\[
\frac{1}{t_{n}}(\beta_{n}-\overline{y})\rightarrow v.
\]
Taking into account the convexity of $\operatorname*{cone}M,$ for every $n,$
\[
\beta_{n}-\overline{y}=t_{n}v_{n}+\tau_{n}w_{n}\in\operatorname*{cone}%
M+\operatorname*{cone}M=\operatorname*{cone}M.
\]
Summing up,%
\[
\frac{1}{t_{n}}(\beta_{n}-\overline{y})\overset{\operatorname*{cone}%
M}{\longrightarrow}v,
\]
whence $v\in T_{B}^{M}(F\left(  \overline{x}\right)  ,\overline{y}).$%
\hfill$\square$

\bigskip

Consider now the situation when $G:X\rightrightarrows Z$ is a set-valued map,
$Q\subset Z$ is a closed convex and pointed cone and the set of restrictions
for $(P)$ is $A:=\{x\in X\mid0\in G(x)+Q\}.$ This is a standard situation
which encompasses the classical case where one has equalities and inequalities
constraints. The following result holds.

\begin{pr}
\label{prop_calc_t}Let $\emptyset\neq L\subset S_{X},$ $\emptyset\neq N\subset
S_{Z}\ $be closed sets, take $\overline{x}\in A$ (meaning that there is
$\overline{z}\in G(\overline{x})\cap-Q$), and define the set-valued map
$\mathcal{E}_{G}:X\rightrightarrows Z,$ $\mathcal{E}_{G}(x)=G(x)+Q.$ Suppose
that $\mathcal{E}_{G}$ is directionally metric subregular at $(\overline
{x},0)$ with respect to $L$ and $N.$ If $\operatorname*{cone}L$ is convex then
$u\in T_{B}^{L}(A,\overline{x})$ iff $0\in D_{B}^{L,N}\mathcal{E}%
_{G}(\overline{x},0)(u).$ Moreover, if $Q\cap S_{Z}\subset N$ and
$\operatorname*{cone}N$ is convex then for every $u\in X,$
\[
D_{B}^{L,N}G(\overline{x},\overline{z})(u)+Q\subset D_{B}^{L,N}\mathcal{E}%
_{G}(\overline{x},0)(u).
\]

\end{pr}

\noindent\textbf{Proof. }We remark that $A=\mathcal{E}_{G}^{-1}(0)$, whence,
by Propositions \ref{prop_echiv_mr_c} and \ref{prop_calm},%
\[
T_{B}^{L}(A,\overline{x})=D_{B}^{N,L}\mathcal{E}_{G}^{-1}(0,\overline{x})(0),
\]
whence $u\in T_{B}^{L}(A,\overline{x})$ iff $u\in D_{B}^{N,L}\mathcal{E}%
_{G}^{-1}(0,\overline{x})(0)$ iff $0\in D_{B}^{L,N}\mathcal{E}_{G}%
(\overline{x},0)(u).$

Now, for the second part, take $w\in D_{B}^{L,N}G(\overline{x},\overline
{z})(u)+Q.$ Then there exist $q\in Q$ and $(u_{n})\overset
{\operatorname*{cone}L}{\longrightarrow}u,(w_{n})\overset{\operatorname*{cone}%
N}{\longrightarrow}w-q,(t_{n})\overset{(0,\infty)}{\longrightarrow}0$ such
that for all $n,$
\[
\overline{z}+t_{n}w_{n}\in G(\overline{x}+t_{n}u_{n}),
\]
whence%
\[
t_{n}(w_{n}+q)\in G(\overline{x}+t_{n}u_{n})-\overline{z}+t_{n}q\subset
\mathcal{E}_{G}(\overline{x}+t_{n}u_{n}).
\]
But, $w_{n}+q\rightarrow w$ and for every $n,$ $w_{n}+q\in\operatorname*{cone}%
N+Q\subset\operatorname*{cone}N,$ whence $w\in D_{B}^{L,N}\mathcal{E}%
_{G}(\overline{x},0)(u)$.\hfill$\square$

\begin{pr}
\label{prop_int_vid}Suppose that $\operatorname*{int}K\neq\emptyset$ and
$(\overline{x},\overline{y})\in\operatorname*{Gr}F$ is a local weak
directional Pareto minimum point for $F$ on $A:=\mathcal{E}_{G}^{-1}(0)$ with
respect to a closed nonempty set $L\subset S_{X}.$ Consider $\overline{z}\in
G(\overline{x})\cap-Q$ and $\emptyset\neq N\subset S_{Z}\ $a closed set$.$
Moreover, suppose that $Q\cap S_{Z}\subset N,$ $\operatorname*{cone}L$ and
$\operatorname*{cone}N$ are convex, and $\mathcal{E}_{G}$ is directionally
metric subregular at $(\overline{x},0)$ with respect to $L$ and $N.$ Then%
\[
\left\{  (v,w)\mid\exists u\in X,v\in D_{D}F(\overline{x},\overline
{y})(u),w\in D_{B}^{L,N}G(\overline{x},\overline{z})(u)\right\}
\cap(-\operatorname*{int}K\times-Q)=\emptyset.
\]

\end{pr}

\noindent\textbf{Proof. }The result follows by using successively Propositions
\ref{prop_weak} and \ref{prop_calc_t}.\hfill$\square$

\bigskip

Let us to specialize, in two steps, the ideas above to the classical smooth
case of optimization problems with single-valued maps. First, suppose that
$F:=f$ and $G:=g$ are continuously Fr\'{e}chet differentiable functions. Then
taking a point $\overline{x}\in A=\mathcal{E}_{G}^{-1}(0)$ it is easy to see
that for all $u\in X,$ $D_{D}f(\overline{x})(u)=\left\{  \nabla f(\overline
{x})(u)\right\}  ,$ while
\[
D_{B}^{L,S_{Z}}g(\overline{x})(u)=\left\{
\begin{array}
[c]{l}%
\left\{  \nabla g(\overline{x})(u)\right\}  ,\text{ if }u\in
\operatorname*{cone}L\\
\emptyset,\text{ if }u\notin\operatorname*{cone}L.
\end{array}
\right.
\]
Then we get the following Fritz John and Karush-Kuhn-Tucker type result.

\begin{thm}
\label{th_nec_1}Suppose that $\operatorname*{int}K\neq\emptyset$ and
$\overline{x}\in A:=\mathcal{E}_{g}^{-1}(0)$ is a local weak directional
Pareto minimum point for $f$ on $A$ with respect to $L.$ Moreover, suppose
that $\operatorname*{cone}L$ is convex, and $\mathcal{E}_{g}$ is directionally
metric subregular at $(\overline{x},0)$ with respect to $L$ and $S_{Z}.$ Then,
in either of the following conditions:

(i) $\operatorname*{int}Q\neq\emptyset$ or $\operatorname*{int}\{(\nabla
f(\overline{x})(u),\nabla g(\overline{x})(u))\mid u\in\operatorname*{cone}%
L\}\neq\emptyset;$

(ii) $Y$ and $Z$ are finite dimensional spaces,

\noindent there exist $y^{\ast}\in K^{+},$ $z^{\ast}\in Q^{+},$ $(y^{\ast
},z^{\ast})\neq0$ such that for every $u\in\operatorname*{cone}L,$%
\[
\left(  y^{\ast}\circ\nabla f(\overline{x})+z^{\ast}\circ\nabla g(\overline
{x})\right)  (u)\geq0.
\]
If, moreover, there exists $u\in\operatorname*{cone}L$ such that $\nabla
g(\overline{x})(u)\in\operatorname*{int}Q\neq\emptyset$ or $\nabla
g(\overline{x})(\operatorname*{cone}L)=Z$ then $y^{\ast}\neq0$.
\end{thm}

\noindent\textbf{Proof.} According to Proposition \ref{prop_int_vid} and the
subsequent discussion,
\[
\{(\nabla f(\overline{x})(u),\nabla g(\overline{x})(u))\mid u\in
\operatorname*{cone}L\}\cap(-\operatorname*{int}K\times-Q)=\emptyset.
\]
Notice that $\{(\nabla f(\overline{x})(u),\nabla g(\overline{x})(u))\mid
u\in\operatorname*{cone}L\}$ is a convex set and both (i) and (ii) ensure the
possibility to apply a separation result for convex sets. Therefore, there
exist $y^{\ast}\in Y^{\ast},$ $z^{\ast}\in Z^{\ast},$ $(y^{\ast},z^{\ast}%
)\neq0$ such that for every $u\in\operatorname*{cone}L,$ $k\in
\operatorname*{int}K,$ $q\in Q,$ one has
\[
\left(  y^{\ast}\circ\nabla f(\overline{x})+z^{\ast}\circ\nabla g(\overline
{x})\right)  (u)\geq-y^{\ast}(k)-z^{\ast}(q).
\]
Standard arguments yield $y^{\ast}\in K^{+},$ $z^{\ast}\in Q^{+}$ and
\[
\left(  y^{\ast}\circ\nabla f(\overline{x})+z^{\ast}\circ\nabla g(\overline
{x})\right)  (u)\geq0.
\]
for every $u\in\operatorname*{cone}L$.

If one supposes that $y^{\ast}=0$ then the relation above and the either of
the final assumptions give $z^{\ast}=0$, which contradicts $(y^{\ast},z^{\ast
})\neq0.$\hfill$\square$

\bigskip

A similar but different result could be done taking into account the special
structure of this case, using directly Proposition \ref{prop_weak}, and some
results one can find in literature concerning the calculus of Bouligand
tangent cone to the counter image of a set through a differentiable mapping.
Let us recall some facts from \cite{DSJOGO2013}. Let $f:X\rightarrow Y$ be a
function and $D\subset X$ be a nonempty closed set. One says that $f$ is
metrically subregular at $(\overline{x},f(\overline{x}))\in D\times Y$ with
respect to $D$ if there exist $s>0,$ $\mu>0$ s.t. for every $u\in
B(\overline{x},s)\cap D$
\[
d(u,f^{-1}(f(\overline{x}))\cap D)\leq\mu\left\Vert f(\overline{x}%
)-f(u)\right\Vert .
\]

In fact, the above notion coincides with that of calmness of the set-valued
map $y\rightrightarrows f^{-1}(y)\cap D$ at $(f(\overline{x}),\overline{x})$
(see, for instance, \cite[Section 3H]{DontRock2009b}). One of the main results
in \cite{DSJOGO2013} reads as follows.

\begin{thm}
\label{th_teor}Let $X,Y$ be Banach spaces, $D\subset X,E\subset Y$ be closed
sets, $\varphi:X\rightarrow Y$ be a continuously Fr\'{e}chet differentiable
map and $\overline{x}\in D\cap\varphi^{-1}(E).$ Suppose that $\psi:X\times
Y\rightarrow Y,$ $\psi(x,y):=\varphi(x)-y$ is metrically subregular at
$(\overline{x},\varphi(\overline{x}),0)$ with respect to $D\times E.$ Then%
\[
T_{U}(D,\overline{x})\cap\nabla\varphi(\overline{x})^{-1}(T_{B}(E,\varphi
(\overline{x})))\subset T_{B}(D\cap\varphi^{-1}(E),\overline{x}),
\]
where $T_{U}(D,\overline{x})$ denotes the Ursescu tangent cone to $D$ at
$\overline{x}$, that is,%
\[
T_{U}(D,\overline{x}):=\left\{  u\in X\mid\forall(t_{n})\overset{(0,\infty
)}{\longrightarrow}0,\exists(u_{n})\rightarrow u\text{ such that for all
}n,\text{ }\overline{x}+t_{n}u_{n}\in D\right\}  .
\]

\end{thm}

\bigskip

Coming back to our case, we have $X:=X,$ $Y:=Z,$ $D:=\overline{x}%
+\operatorname*{cone}L,$ $E:=-Q,$ $\varphi:=g.$ We have seen that $T_{B}%
^{L}(g^{-1}(-Q),\overline{x})=T_{B}((\overline{x}+\operatorname*{cone}L)\cap
g^{-1}(-Q)),\overline{x}).$ With these identifications, we get the next result.

\begin{thm}
\label{T3.11}Suppose that $X,Z$ are Banach spaces, $\operatorname*{int}%
K\neq\emptyset$ and $\overline{x}\in g^{-1}(-Q)$ is a local weak directional
Pareto minimum point for $f$ on $g^{-1}(-Q)$ with respect to $L.$ Moreover,
suppose that $\psi:X\times Z\rightarrow Z,$ $\psi(x,z):=g(x)-z$ is metrically
subregular at $(\overline{x},g(\overline{x}),0)$ with respect to
$(\overline{x}+\operatorname*{cone}L)\times-Q.$ Then for all $u\in
\operatorname*{cone}L$ with $\nabla g(\overline{x})(u)\in T_{B}(-Q,g(\overline
{x})),$%
\[
\nabla f(\overline{x})(u)\notin-\operatorname*{int}K
\]

\end{thm}

\noindent\textbf{Proof.} According to Theorem \ref{th_teor},
\[
T_{U}(\overline{x}+\operatorname*{cone}L,\overline{x})\cap\nabla
g(\overline{x})^{-1}(T_{B}(-Q,g(\overline{x})))\subset T_{B}^{L}%
(g^{-1}(-Q),\overline{x}),
\]
whence%
\[
\operatorname*{cone}L\cap\nabla g(\overline{x})^{-1}(T_{B}(-Q,g(\overline
{x})))\subset T_{B}^{L}(g^{-1}(-Q),\overline{x}).
\]
By Proposition \ref{prop_weak},
\[
\nabla f(\overline{x})(u)\notin-\operatorname*{int}K
\]
for all $u\in\operatorname*{cone}L\cap\nabla g(\overline{x})^{-1}%
(T_{B}(-Q,g(\overline{x}))),$ whence the conclusion.\hfill$\square$

\bigskip

Furthermore, we consider the case where $Y=\mathbb{R}^{k}$ $(k\geq1),$
$Z=\mathbb{R}^{p}$ $(p\geq1),$ $Q=\mathbb{R}_{+}^{m}\times\{0\}^{n}$ with
$m+n=p,$ and $f,g$ are Fr\'{e}chet differentiable. This means that we are
dealing with a vectorial optimization problem with finitely many inequalities
and equalities constraints. Let us denote by $\mu_{i}$ with $i\in
\overline{1,m}$ the first $m$ coordinates functions of $g$ and by $\nu_{j}$
with $j\in\overline{1,n}$ the next $n$ coordinates functions of $g.$

For the next step of our approach, we use the Gerstewitz functional in the
special case when the ordering cone has nonempty interior. The next result
combines \cite[Theorem 2.3.1]{GRTZ-S} and \cite[Lemma 2.1]{D-PJO}.

\begin{thm}
\label{Lemma_funct}Let $K\subset Y$ be a closed convex cone with nonempty
interior. Then for every $e\in\operatorname{int}K$ the functional
$s_{K,e}:Y\rightarrow\mathbb{R}$ given by
\begin{equation}
s_{K,e}(y)=\inf\{\lambda\in\mathbb{R}\mid\lambda e\in y+K\} \label{sep_subln}%
\end{equation}
is convex continuous and for every $\lambda\in\mathbb{R},$
\begin{equation}
\{y\in Y\mid s_{K,e}(y)<\lambda\}=\lambda e-\operatorname*{int}K,\text{ and
}\{y\in Y\mid s_{K,e}(y)=\lambda\}=\lambda e-\operatorname{bd}K.
\label{level_str}%
\end{equation}

Moreover, $s_{K,e}$ is sublinear, $K-$monotone, and for every $u\in Y,$ the
Fenchel (convex) subdifferential $\partial s_{K,e}(u)$ is nonempty and
\begin{equation}
\partial s_{K,e}(u)=\{v^{\ast}\in K^{+}\mid v^{\ast}(e)=1,v^{\ast}%
(u)=s_{K,e}(u)\}. \label{subd_e}%
\end{equation}

\end{thm}

In this notation we have the next result.

\begin{thm}
Suppose that $X$ is a Banach space, $\operatorname*{int}K\neq\emptyset$ and
$\overline{x}\in g^{-1}(-Q)$ is a local weak directional Pareto minimum point
for $f$ on $g^{-1}(-Q)$ with respect to $L.$ Suppose that:

(i) $\operatorname*{cone}L$ is convex;

(ii) $\psi:X\times Z\rightarrow Z,$ $\psi(x,z):=g(x)-z$ is metrically
subregular at $(\overline{x},g(\overline{x}),0)$ with respect to
$(\overline{x}+\operatorname*{cone}L)\times-Q;$

(iii) $\nabla\nu\left(  \overline{x}\right)  \left(  X\right)  =\mathbb{R}%
^{n},$ where $\nu:=\left(  \nu_{1},\nu_{2},...,\nu_{n}\right)  ;$

(iv) there exists $\overline{u}\in\operatorname*{int}\operatorname*{cone}L$
such that $\nabla\mu_{i}\left(  \overline{x}\right)  \left(  \overline
{u}\right)  <0$ for any $i\in I(\overline{x}):=\{i\in\overline{1,m}\mid\mu
_{i}(\overline{x})=0\}$ and $\nabla\nu\left(  \overline{x}\right)  \left(
\overline{u}\right)  =0.$

Then there exist $y^{\ast}\in K^{+}\mathbb{\setminus\{}0\mathbb{\}},$
$\lambda_{i}\geq0$ for $i\in\overline{1,m}$ and $\tau_{j}\in\mathbb{R}$ for
$j\in\overline{1,n}$ such that%
\begin{equation}
0\in y^{\ast}\circ\nabla f(\overline{x})+\sum\limits_{i=1}^{m}\lambda
_{i}\nabla\mu_{i}(\overline{x})+\sum\limits_{j=1}^{n}\tau_{j}\nabla\nu
_{j}(\overline{x})+L^{-} \label{kkt1}%
\end{equation}
and
\begin{equation}
\lambda_{i}\mu_{i}(\overline{x})=0,\forall i\in\overline{1,m}. \label{kkt2}%
\end{equation}

\end{thm}

\noindent\textbf{Proof.} Clearly, in this case $u\in\nabla g(\overline
{x})^{-1}(T_{B}(-Q,g(\overline{x})))$ amounts to say that $\nabla\mu
_{i}(\overline{x})(u)\leq0$ for any $i\in I(\overline{x})$ and $\nabla\nu
_{j}(\overline{x})(u)=0$ for any $j\in\overline{1,n}$.

Using Theorem \ref{T3.11} (all its assumptions hold) we get that
\[
\nabla f(\overline{x})(u)\notin-\operatorname*{int}K
\]
for all $u\in\operatorname*{cone}L$ with $\nabla\mu_{i}(\overline{x})(u)\leq0$
for any $i\in I(\overline{x}),$ and $\nabla\nu_{j}(\overline{x})(u)=0$ for any
$j\in\overline{1,n}.$

We conclude that $s_{K,e}\left(  \nabla f(\overline{x})(u)\right)  \geq0$ for
all $u$ satisfying the above conditions and this means that $u=0$ is a minimum
point for the scalar problem
\[
\min s_{K,e}\left(  \nabla f(\overline{x})(u)\right)  \text{ s.t. }%
u\in\operatorname*{cone}L,\nabla\mu_{i}(\overline{x})(u)\leq0,\forall i\in
I(\overline{x}),\nabla\nu_{j}(\overline{x})(u)=0,\forall j\in\overline{1,n}.
\]
Since $\operatorname*{cone}L$ is convex, this is a convex problem, whence,
from \cite[Theorem 2.9.6]{Zal2002}, there exist $\lambda_{i}\geq0$ for $i\in
I(\overline{x})$ and $\tau_{j}\in\mathbb{R}$ for $j\in\overline{1,n}$ such
that
\[
0\in\partial\left(  s_{K,e}\circ\nabla f(\overline{x})+\iota
_{\operatorname*{cone}L}+\sum_{i\in I(\overline{x})}\lambda_{i}\nabla\mu
_{i}(\overline{x})+\sum\limits_{j=1}^{n}\tau_{j}\nabla\nu_{j}(\overline
{x})\right)  (0),
\]
where $\iota$ denotes the indicator function. Finally, using (\ref{subd_e}),
and taking $\lambda_{i}:=0$ for $i\in\overline{1,m}\mathbb{\setminus
}I(\overline{x}),$ we get the existence of $y^{\ast}\in K^{+}\mathbb{\setminus
\{}0\mathbb{\}}$ such that%
\[
0\in y^{\ast}\circ\nabla f(\overline{x})+\sum\limits_{i=1}^{m}\lambda
_{i}\nabla\mu_{i}(\overline{x})+\sum\limits_{j=1}^{n}\tau_{j}\nabla\nu
_{j}(\overline{x})+L^{-}%
\]
and
\[
\lambda_{i}\mu_{i}(\overline{x})=0,\forall i\in\overline{1,m},
\]
whence the conclusion.\hfill$\square$

\begin{rmk}
Observe that in the simplest case of a derivable real-valued function
$f:\mathbb{R\rightarrow R},$ if $\overline{x}$ is a directional minimum with
respect to $L:=\{+1\}$ (without constraints) the above theorem reduces to
$-f(\overline{x})\in L^{-}$ which is exactly $f^{\prime}(\overline{x})\geq0,$
as discussed before.
\end{rmk}

Our aim now is to derive sufficient conditions for a point $\overline{x}\in
g^{-1}\left(  -Q\right)  $ to be a local weak directional Pareto minimum
point. In order to formulate such conditions we use, besides the convexity
notion for scalar functions, a generalized convexity concept. Namely, we use
the following well-known concept: one says that $F:X\rightrightarrows Y$ is
$K-$convex if for any $\lambda\in\left(  0,1\right)  ,$ and any $x,y\in X,$
one has%
\[
\lambda F(x)+(1-\lambda)F(y)\subset F(\lambda x+(1-\lambda)y)+K.
\]

\begin{pr}
Suppose that $X$ is a Banach space, $\operatorname*{int}K\neq\emptyset,$
$\operatorname*{cone}L$ is convex, $f$ is $K-$convex, $\mu_{i},i\in
\overline{1,m}$, are convex and $\nu_{j},j\in\overline{1,n},$ are affine. If
there exist $\left(  \lambda,\tau\right)  \in\mathbb{R}_{+}^{m}\times
\mathbb{R}^{n}$ and $y^{\ast}\in K^{+}\setminus\left\{  0\right\}  $ such that
(\ref{kkt1}) and (\ref{kkt2}) hold, then $\overline{x}$ is a global weak
directional Pareto minimum point for $f$ on $g^{-1}\left(  -Q\right)  $ with
respect to $L.$
\end{pr}

\noindent\textbf{Proof.} By relation (\ref{kkt1}), we immediately get that%
\begin{align*}
0 &  \in\nabla\left(  y^{\ast}\circ f+\sum_{i=1}^{m}\lambda_{i}\mu_{i}%
+\sum_{j=1}^{n}\tau_{j}\nu_{j}\right)  \left(  \overline{x}\right)  +N\left(
\operatorname*{cone}L,0\right)  \\
&  =\nabla\left(  y^{\ast}\circ f+\sum_{i=1}^{m}\lambda_{i}\mu_{i}+\sum
_{j=1}^{n}\tau_{j}\nu_{j}\right)  \left(  \overline{x}\right)  +N\left(
\overline{x}+\operatorname*{cone}L,\overline{x}\right)  .
\end{align*}
Consider the convex optimization problem%
\begin{equation}
\min\left(  \left(  y^{\ast}\circ f\right)  \left(  x\right)  +\sum_{i=1}%
^{m}\lambda_{i}\mu_{i}\left(  x\right)  +\sum_{j=1}^{n}\tau_{j}\nu_{j}\left(
x\right)  \right)  ,\quad x\in\overline{x}+\operatorname*{cone}L.\label{kkt3}%
\end{equation}
We hence obtain, by virtue of \cite[Theorem 2.9.1]{Zal2002}, that
$\overline{x}$ is a global minimum point for the above problem. Note that, for
all feasible points $x\in g^{-1}\left(  -Q\right)  ,$ we have%
\[
\sum_{i=1}^{m}\lambda_{i}\mu_{i}\left(  x\right)  +\sum_{j=1}^{n}\tau_{j}%
\nu_{j}\left(  x\right)  =\sum_{i=1}^{m}\lambda_{i}\mu_{i}\left(  x\right)
\leq0.
\]
Using (\ref{kkt2}), it follows that, given any $x\in\left(  \overline
{x}+\operatorname*{cone}L\right)  \cap g^{-1}\left(  -Q\right)  ,$%
\[
\left(  y^{\ast}\circ f\right)  \left(  x\right)  \geq\left(  y^{\ast}\circ
f\right)  \left(  \overline{x}\right)  ,
\]
that is%
\[
y^{\ast}\left(  f\left(  x\right)  -f\left(  \overline{x}\right)  \right)
\geq0.
\]
Now, since $y^{\ast}\in K^{+}\mathbb{\setminus}\left\{  0\right\}  ,$ the
inequality above gives $f\left(  x\right)  -f\left(  \overline{x}\right)
\notin-\operatorname*{int}K,$ i.e., the conclusion.\hfill$\square$

\subsection{Optimality conditions using normal cones}

In order to tackle the question of optimality conditions for directional
minima in terms of generalized differentiation objects in dual spaces, we
recall some notions and results concerning Fr\'{e}chet and limiting
(Mordukhovich) generalized differentiation (see \cite{Mor-2006-v1} for details).

Consider $S$ a nonempty subset$\ $of a Banach space $X$ and $x\in S.$ Then for
every $\varepsilon\geq0,$ the set of $\varepsilon-$normals to $S$ at $x$ is
defined by%
\[
\widehat{N}_{\varepsilon}(S,x)=\left\{  x^{\ast}\in X^{\ast}\mid
\underset{u\overset{S}{\rightarrow}x}{\lim\sup}\frac{x^{\ast}(u-x)}{\left\Vert
u-x\right\Vert }\leq\varepsilon\right\}  ,
\]
where $u\overset{S}{\rightarrow}x$ means that $u\rightarrow x$ and $u\in S.$
The set $\widehat{N}_{0}(S,x)$ is denoted by $\widehat{N}(S,x)$ and it is
called the Fr\'{e}chet normal cone to $S$ at $x.$

Let $\overline{x}\in S.$ The Mordukhovich normal cone to $S$ at $\overline{x}$
is given by
\[
N(S,\overline{x})=\{x^{\ast}\in X^{\ast}\mid\exists\varepsilon_{n}%
\overset{\left(  0,\infty\right)  }{\longrightarrow}0,x_{n}\overset
{S}{\rightarrow}\overline{x},x_{n}^{\ast}\overset{w^{\ast}}{\rightarrow
}x^{\ast},x_{n}^{\ast}\in\widehat{N}_{\varepsilon_{n}}(S,x_{n}),\forall
n\in\mathbb{N}\}.
\]
Up to the end of this section, we consider that all the involved spaces are
Asplund, unless otherwise stated. In this context, if $S\subset X$ is closed
around $\overline{x}$, the formula for the Mordukhovich normal cone takes the
following form:%
\[
N(S,\overline{x})=\{x^{\ast}\in X^{\ast}\mid\exists x_{n}\overset
{S}{\rightarrow}\overline{x},x_{n}^{\ast}\overset{w^{\ast}}{\rightarrow
}x^{\ast},x_{n}^{\ast}\in\widehat{N}(S,x_{n}),\forall n\in\mathbb{N}\}.
\]

For the set-valued map\ $F:X\rightrightarrows Y,$ its Fr\'{e}chet coderivative
at $(\overline{x},\overline{y})\in\operatorname*{Gr}F$ is the set-valued map
$\widehat{D}^{\ast}F(\overline{x},\overline{y}):Y^{\ast}\rightrightarrows
X^{\ast}$ given by
\[
\widehat{D}^{\ast}F(\overline{x},\overline{y})(y^{\ast})=\{x^{\ast}\in
X^{\ast}\mid(x^{\ast},-y^{\ast})\in\widehat{N}(\operatorname{Gr}%
F,(\overline{x},\overline{y}))\}.
\]
In the same way, the Mordukhovich coderivative of $F$ at $(\overline
{x},\overline{y})$ is the set-valued map $D^{\ast}F(\overline{x},\overline
{y}):Y^{\ast}\rightrightarrows X^{\ast}$ given by
\[
D^{\ast}F(\overline{x},\overline{y})(y^{\ast})=\{x^{\ast}\in X^{\ast}%
\mid(x^{\ast},-y^{\ast})\in N(\operatorname{Gr}F,(\overline{x},\overline
{y}))\}.
\]
As usual, when $F=f$ is a function, since $\overline{y}\in F\left(
\overline{x}\right)  $ means $\overline{y}=f\left(  \overline{x}\right)  ,$ we
write $\widehat{D}^{\ast}f\left(  \overline{x}\right)  $ for $\widehat
{D}^{\ast}f\left(  \overline{x},\overline{y}\right)  ,$ and similarly for
$D^{\ast}.$

Notice that for a convex set $S\subset X$ one has that%
\[
N(S,\overline{x})=\{x^{\ast}\in X^{\ast}\mid x^{\ast}(x-\overline{x}%
)\leq0,\forall x\in S\}
\]
and this cone coincides with the negative polar of $T_{B}(S,\overline{x}).$

If $S\subset X$ is closed around $\overline{x}\in S,$ one says that $S$ is
sequentially normally compact (SNC, for short) at $\overline{x}$ if
\[
\left[  x_{n}\overset{S}{\rightarrow}\overline{x},x_{n}^{\ast}\overset
{w^{\ast}}{\rightarrow}0,x_{n}^{\ast}\in\widehat{N}(S,x_{n})\right]
\Rightarrow x_{n}^{\ast}\rightarrow0.
\]
In the case where $S=C$ is a closed convex cone, the (SNC) property at $0$ is
equivalent to%
\[
\left[  (x_{n}^{\ast})\subset C^{+},x_{n}^{\ast}\overset{w^{\ast}}%
{\rightarrow}0\right]  \Rightarrow x_{n}^{\ast}\rightarrow0.
\]
In particular, if $\operatorname*{int}C\neq\emptyset,$ then $C$ is (SNC) at
$0.$

Let $f:X\rightarrow\mathbb{R\cup\{+\infty\}}$ be finite at $\overline{x}\in
X\ $and lower semicontinuous around $\overline{x};$ the Fr\'{e}chet
subdifferential of $f$ at $\overline{x}$ is defined by
\[
\widehat{\partial}f(\overline{x})=\{x^{\ast}\in X^{\ast}\mid(x^{\ast}%
,-1)\in\widehat{N}(\operatorname*{epi}f,(\overline{x},f(\overline{x})))\},
\]
where $\operatorname*{epi}f$ denotes the epigraph of $f.$ The Mordukhovich
subdifferential of $f$ at $\overline{x}$ is given by%
\[
\partial f(\overline{x})=\{x^{\ast}\in X^{\ast}\mid(x^{\ast},-1)\in
N(\operatorname*{epi}f,(\overline{x},f(\overline{x})))\}.
\]
It is well-known that if $f$ is a convex function, then $\widehat{\partial
}f(\overline{x})$ and $\partial f(\overline{x})$ coincide with the Fenchel
subdifferential. However, in general, $\widehat{\partial}f(\overline
{x})\subset\partial f(\overline{x}),$ and the following generalized Fermat
rule holds: if $\overline{x}\in X$ with $f(\overline{x})<+\infty$ is a local
minimum point for $f:X\rightarrow\mathbb{R\cup\{+\infty\}}$, then
$0\in\widehat{\partial}f(\overline{x}).$

Consider now some subsets $C_{1},...,C_{k}$ of $X$ ($k\in\mathbb{N\setminus
}\{0,1\}$)$.$ Take $\overline{x}\in C_{1}\cap...\cap C_{k}$ and suppose that
all the sets $C_{i}$, $i\in\overline{1,k}$ are closed around $\overline{x}.$
One says that $C_{1},...,C_{k}$ are allied at $\overline{x}$ if for every
$(x_{in})\overset{C_{i}}{\rightarrow}\overline{x},$ $x_{in}^{\ast}\in
\widehat{N}(C_{i},x_{in}),$ $i\in\overline{1,k},$ the relation $(x_{1n}^{\ast
}+...+x_{kn}^{\ast})\rightarrow0$ implies $(x_{in}^{\ast})\rightarrow0$ for
every $i\in\overline{1,k}$. The concept of alliedness was introduced by Penot
and his coauthors in \cite{penot-allied} and \cite{penot+altii-aliere} in
order to get a calculus rule for the Fr\'{e}chet normal cone to the
intersection of sets. More precisely, if the subsets $C_{1},...,C_{k}$ are
allied at $\overline{x},$ then there exists $r>0$ such that, for every
$\varepsilon>0$ and every $x\in\lbrack C_{1}\cap...\cap C_{k}]\cap
B_{X}(\overline{x},r)$, there exist $x_{i}\in C_{i}\cap B_{X}(x,\varepsilon)$,
$i\in\overline{1,k}$ such that%
\[
\widehat{N}(C_{1}\cap...\cap C_{k},x)\subset\widehat{N}(C_{1},x_{1}%
)+...+\widehat{N}(C_{k},x_{k})+\varepsilon D_{X^{\ast}}.
\]

In what follows we use the results concerning the theory of generalized
differentiation built on these objects directly at the places we need them,
without separate quotation.

\bigskip

We discuss next a concept of directional openness at the reference point of a
certain multifunction. We recall that the classical concept of openness proven
to be useful for the announced aim by means of the incompatibility between
this property and the Pareto minimality (see, e.g., \cite{DS-Opt} for details).

In fact, the directional openness we consider here is related to several other
notions introduced in \cite{DPS}, and to the concept of directional calmness
already used in the previous subsection.

Consider a multifunction $F:X\rightrightarrows Y,$ a point $(\overline
{x},\overline{y})\in\operatorname*{Gr}F,$ and $\emptyset\neq L\subset S_{X},$
$\emptyset\neq M\subset S_{Y}.$ One says $F$ is directionally open at
$(\overline{x},\overline{y})$ with respect to $L$ and $M$ if for any
$\varepsilon>0,$ there exists $r>0$ such that%
\[
B(\overline{y},r)\cap\lbrack\overline{y}-\operatorname*{cone}M]\subset
F(B(\overline{x},\varepsilon)\cap\lbrack\overline{x}+\operatorname*{cone}L]).
\]

When $F$ is single-valued, for simplicity, we sometimes omit $\overline{y}$ in
the definition above and we say that $F$ is directionally open at
$\overline{x},$ instead of directionally open at $(\overline{x},f(\overline
{x}))$.

\begin{pr}
\label{pr_not_dir_op}If $(\overline{x},\overline{y})\in\operatorname*{Gr}F$ is
a local directional Pareto minimum point for $F$ with respect to $L,$ then for
every $C\subset S_{Y}$ with $C\cap\left(  K\setminus-K\right)  \neq\emptyset,$
the set-valued map $\mathcal{E}_{F}:X\rightrightarrows Y,$ given by
$\mathcal{E}_{F}(x):=F(x)+K$ is not directionally open at $(\overline
{x},\overline{y})$ with respect to $L$ and $C.$ In particular, $F$ is not
directionally open at $(\overline{x},\overline{y})$ with respect to $L$ and
$C.$
\end{pr}

\noindent\textbf{Proof. }Suppose, by contradiction, that for $\varepsilon>0$
involved in the definition of the minimality of $(\overline{x},\overline{y}),$
there exists $r>0$ such that%
\[
B(\overline{y},r)\cap\lbrack\overline{y}-\operatorname*{cone}C]\subset
\mathcal{E}_{F}(B(\overline{x},\varepsilon)\cap\lbrack\overline{x}%
+\operatorname*{cone}L]).
\]
By subtracting $\overline{y}$ on both sides, according to the hypothesis, one
has that%
\begin{align*}
\left[  B(0,r)\cap-C\right]  \cap-K  &  \subset\left[  \mathcal{E}%
_{F}(B(\overline{x},\varepsilon)\cap\lbrack\overline{x}+\operatorname*{cone}%
L])-\overline{y}\right]  \cap-K\\
&  =\left[  F(B(\overline{x},\varepsilon)\cap\lbrack\overline{x}%
+\operatorname*{cone}L])+K-\overline{y}\right]  \cap-K\subset K.
\end{align*}
Passing to the conic hull, this yields%
\[
-\operatorname*{cone}C\cap-K\subset K,
\]
which contradicts the fact that $C\cap\left(  K\setminus-K\right)
\neq\emptyset$. So $\mathcal{E}_{F}$ is not directionally open at
$(\overline{x},\overline{y})$ with respect to $L$ and $C.$ Since
$F(x)\subset\mathcal{E}_{F}(x)$ for any $x,$ the same conclusion holds for $F$
as well.\hfill$\square$

\bigskip

Before obtaining necessary optimality conditions, we remark that a converse of
Proposition \ref{pr_not_dir_op} can be done if one considers a (generalized)
convex framework.

\begin{pr}
Suppose that $F$ is $K-$convex and for every $u\in K\cap S_{Y},$
$\mathcal{E}_{F}$ is not directionally open with respect to $L$ and
$M:=\left\{  u\right\}  $ at $(\overline{x},\overline{y})\in\operatorname*{Gr}%
F.$ Then $(\overline{x},\overline{y})$ is a local directional Pareto minimum
point of $F$ with respect to $L.$
\end{pr}

\noindent\textbf{Proof.} Suppose, by contradiction, that $(\overline
{x},\overline{y})$ is not a local directional Pareto minimum point of $F$ with
respect to $L.$ Then for every $r>0,$ there is $y_{r}\in F(B(\overline
{x},r)\cap\lbrack\overline{x}+\operatorname*{cone}L])\cap\left[  \overline
{y}-K\right]  $ such that $y_{r}\notin\overline{y}+K.$ Denote $\overline
{k}:=\overline{y}-y_{r}\in K\setminus-K$ and consider $x_{r}\in B(\overline
{x},r)\cap\left(  \overline{x}+\operatorname*{cone}L\right)  $ such that
$y_{r}\in F(x_{r}).$

Moreover, since $\mathcal{E}_{F}$ is not directionally open with respect to
$L$ and $\left\{  \overline{k}\right\}  $ at $(\overline{x},\overline{y}),$ it
follows that there is $r>0$ such that, for every $\varepsilon>0\ $small
enough, there is $y_{\varepsilon}\in B(\overline{y},\varepsilon)\cap\left[
\overline{y}-\operatorname*{cone}\overline{k}\right]  \subset\left[
\overline{y},y_{r}\right]  $ such that $y_{\varepsilon}\notin\mathcal{E}%
_{F}(B(\overline{x},r)\cap\lbrack\overline{x}+\operatorname*{cone}L])$ (hence,
in particular, $y_{\varepsilon}\neq\overline{y}$ and $y_{\varepsilon}\neq
y_{r}$)$.$

Then, there is $\lambda\in(0,1)$ such that%
\begin{align*}
y_{\varepsilon}  &  =\lambda\overline{y}+(1-\lambda)y_{r}\in\lambda
F(\overline{x})+(1-\lambda)F(x_{r})\subset F(\lambda\overline{x}%
+(1-\lambda)x_{r})+K\\
&  =\mathcal{E}_{F}\left(  \lambda\overline{x}+(1-\lambda)x_{r}\right)
=\mathcal{E}_{F}\left(  \overline{x}+(1-\lambda)(x_{r}-\overline{x})\right)
\subset\mathcal{E}_{F}(B(\overline{x},r)\cap\lbrack\overline{x}%
+\operatorname*{cone}L]),
\end{align*}
a contradiction.\hfill$\square$

\bigskip

Now, we use Proposition \ref{pr_not_dir_op} to get optimality conditions.

\begin{thm}
Suppose that $X$ and $Y\ $are finite dimensional spaces, $(\overline
{x},\overline{y})\in\operatorname*{Gr}F$ is a local directional Pareto minimum
point for $F$ with respect to $L$, $\operatorname*{cone}L$ is convex,
$u\in\operatorname*{int}K\cap S_{Y},$ and the set-valued map $\mathcal{E}%
_{F}:X\rightrightarrows Y\ $has closed graph and is Lipschitz-like around
$(\overline{x},\overline{y})$. Then there exist $x^{\ast}\in X^{\ast},y^{\ast
}\in K^{+}$ with $x^{\ast}(\ell)\geq0$ for all $\ell\in L,$ $y^{\ast}(u)=1$
and
\[
x^{\ast}\in D^{\ast}\mathcal{E}_{F}(\overline{x},\overline{y})(y^{\ast}).
\]

\end{thm}

\noindent\textbf{Proof. }According to Proposition \ref{pr_not_dir_op},
$\mathcal{E}_{F}$ is not directionally open at $(\overline{x},\overline{y})$
with respect to $L$ and $\{u\}.$ Therefore, this is not directionally open
around $(\overline{x},\overline{y})$ with respect to $L$ and $\{u\}$ (in the
sense of \cite[Definition 2.2]{DPS}) and, therefore, the sufficient condition
for directional openness from \cite[Theorem 4.3]{DPS} does not hold. This
means that for all natural numbers $n\neq0,$ there exist $x_{n}^{\ast}\in
X^{\ast},y_{n}^{\ast}\in Y^{\ast},$ $(x_{n},y_{n})\overset{\operatorname*{Gr}%
F}{\longrightarrow}(\overline{x},\overline{y})$ such that $y_{n}^{\ast}(u)=1,$
$n^{-1}>-x_{n}^{\ast}(\ell)$ for all $\ell\in L$ and $x_{n}^{\ast}\in\hat
{D}^{\ast}\mathcal{E}_{F}(x_{n},y_{n})(y_{n}^{\ast}).$ Now \cite[Lemma
3.2]{DS-Opt} ensures that $y_{n}^{\ast}\in K^{+}$ for any $n.$ This, together
with the condition $u\in\operatorname*{int}K$ imply, by using \cite[Lemma
2.2.17]{GRTZ-S}, that the sequence $(y_{n}^{\ast})$ is bounded. The assumed
Lipschitz property of $\mathcal{E}_{F}$ ensures, by means of \cite[Theorem
1.43]{Mor-2006-v1}, that the sequence $(x_{n}^{\ast})$ is bounded too.
Therefore, we can suppose, without loss of generality, that both these
sequences are convergent to some $x^{\ast}\in X^{\ast}$ and $y^{\ast}\in
K^{+},$ respectively. Passing to the limit in the relations satisfied by
$(x_{n}^{\ast})$ and $(y_{n}^{\ast})$ we get, $x^{\ast}(\ell)\geq0$ for all
$\ell\in L,$ $y^{\ast}(u)=1$ and $x^{\ast}\in D^{\ast}\mathcal{E}%
_{F}(\overline{x},\overline{y})(y^{\ast}),$ that is the conclusion.\hfill
$\square$

\begin{rmk}
Observe that, in the case $L=S_{X}$ (that is Pareto minimality) the necessary
optimality condition given by the previous result is the generalized Fermat
rule (see \cite[Theorem 3.11]{DS-Opt}): there exists $y^{\ast}\in
K^{+}\setminus\{0\}$ with
\[
0\in D^{\ast}\mathcal{E}_{F}(\overline{x},\overline{y})(y^{\ast}).
\]

\end{rmk}

We tackle now the case of constrained problems and we have the following result.

\begin{thm}
Let $A\subset X$ and $L\subset S_{X}$ be nonempty closed sets and
$F:X\rightrightarrows Y$ be a set-valued map with $\left(  \overline
{x},\overline{y}\right)  \in\operatorname*{Gr}F\cap\left(  A\times Y\right)  $
such that $\operatorname*{Gr}F$ is closed around $\left(  \overline
{x},\overline{y}\right)  .$ Suppose that the following assertions hold:

(i) $F$ is Lipschitz-like around $\left(  \overline{x},\overline{y}\right)  ;$

(ii) $K\mathbb{\setminus-}K\mathbb{\neq\emptyset}$ and $K$ is (SNC) at $0$;

(iii) the sets $A$ and $\overline{x}+\operatorname*{cone}L$ are allied at
$\overline{x}$.

If $\left(  \overline{x},\overline{y}\right)  $ is a local directional Pareto
minimum point for $F$ on $A$ with respect to the set of directions $L,$ then
there exists $y^{\ast}\in K^{+}\mathbb{\setminus}\left\{  0\right\}  $ such
that%
\[
0\in D^{\ast}F\left(  \overline{x},\overline{y}\right)  \left(  y^{\ast
}\right)  +N\left(  A,\overline{x}\right)  +N\left(  \operatorname*{cone}%
L,0\right)  .
\]

\end{thm}

\noindent\textbf{Proof. }From the hypothesis, there exists a neighborhood
$U\in\mathcal{V}\left(  \overline{x}\right)  $ such that%
\begin{equation}
\left(  F(U\cap A\cap\left(  \overline{x}+\operatorname*{cone}L\right)
)-\overline{y}\right)  \cap-K\subset K \label{ip_min}%
\end{equation}
and there exists $c\in Y$ such that $c\in K\mathbb{\setminus-}K.$ Consider the
following two sets:%
\[
A_{1}=\operatorname*{Gr}F
\]
and%
\[
A_{2}=\left[  \left(  \overline{x}+\operatorname*{cone}L\right)  \cap
A\right]  \times\left(  \overline{y}-K\right)  .
\]
We want to prove that the system $\left\{  A_{1},A_{2},\left(  \overline
{x},\overline{y}\right)  \right\}  $ is an extremal system in $X\times Y$ (see
\cite[Definition 2.1]{Mor-2006-v1}). For this, since the sets $A_{1}$ and
$A_{2}$ are closed around $\left(  \overline{x},\overline{y}\right)  \in
A_{1}\cap A_{2},$ it is sufficient to show the existence of a sequence
$\left(  \left(  x_{n},y_{n}\right)  \right)  _{n}\subset X\times Y$ such that
$\left(  x_{n},y_{n}\right)  \rightarrow\left(  0,0\right)  $ and%
\[
A_{1}\cap\left(  A_{2}-\left(  x_{n},y_{n}\right)  \right)  \cap\left(
U\times Y\right)  =\emptyset,
\]
for all large $n\in\mathbb{N}$. Consider $\left(  x_{n},y_{n}\right)  =\left(
0,\frac{c}{n}\right)  $ with $n\in\mathbb{N\setminus}\left\{  0\right\}  $ and
suppose, by contradiction, that there exist $x\in\left(  \overline
{x}+\operatorname*{cone}L\right)  \cap A\cap U$ and $y\in F\left(  x\right)
\cap\left(  \overline{y}-K-\frac{c}{n}\right)  \subset F\left(  x\right)
\cap\left(  \overline{y}-K\right)  ,$ whence $y-\overline{y}\in\left(
F\left(  x\right)  -\overline{y}\right)  \cap-K$. Now, using (\ref{ip_min}) we
get that $y-\overline{y}\in K$ and since $\overline{y}-y-\frac{c}{n}\in K$ we
arrive at $-c\in K,$ a contradiction. Thus, $\left\{  A_{1},A_{2},\left(
\overline{x},\overline{y}\right)  \right\}  $ is an extremal system in
$X\times Y$ and since $X\times Y$ is an Asplund space we can apply the
approximate extremal principle to this system (see, \cite[Theorem
2.20]{Mor-2006-v1}). Therefore, for every $n\in\mathbb{N\setminus}\left\{
0\right\}  ,$ there exist $\left(  x_{n}^{1},y_{n}^{1}\right)  \in
\operatorname*{Gr}F\cap D\left(  \left(  \overline{x},\overline{y}\right)
,\frac{1}{n}\right)  ,$ $x_{n}^{2}\in\left(  \overline{x}+\operatorname*{cone}%
L\right)  \cap A\cap D\left(  \overline{x},\frac{1}{n}\right)  ,$ $y_{n}%
^{2}\in\left(  \overline{y}-K\right)  \cap D\left(  \overline{y},\frac{1}%
{n}\right)  ,$ $x_{n}^{1\ast}\in X^{\ast},x_{n}^{2\ast}\in X^{\ast}%
,y_{n}^{1\ast}\in Y^{\ast},y_{n}^{2\ast}\in Y^{\ast}$ such that%
\begin{align*}
&  \left(  x_{n}^{1\ast},y_{n}^{1\ast}\right)  \in\widehat{N}\left(
\operatorname*{Gr}F,\left(  x_{n}^{1},y_{n}^{1}\right)  \right)  +\frac{1}%
{n}D_{X^{\ast}\times Y^{\ast}},\\
&  x_{n}^{2\ast}\in\widehat{N}\left(  \left(  \overline{x}%
+\operatorname*{cone}L\right)  \cap A,x_{n}^{2}\right)  +\frac{1}{n}%
D_{X^{\ast}},\\
&  y_{n}^{2\ast}\in\widehat{N}\left(  \overline{y}-K,y_{n}^{2}\right)
+\frac{1}{n}D_{Y^{\ast}}=-\widehat{N}\left(  K,\overline{y}-y_{n}^{2}\right)
+\frac{1}{n}D_{Y^{\ast}}%
\end{align*}
and%
\begin{equation}
x_{n}^{1\ast}+x_{n}^{2\ast}=0,\text{ }y_{n}^{1\ast}+y_{n}^{2\ast}=0,\text{
}\left\Vert \left(  x_{n}^{1\ast},y_{n}^{1\ast}\right)  \right\Vert
+\left\Vert \left(  x_{n}^{2\ast},y_{n}^{2\ast}\right)  \right\Vert =1.
\label{1}%
\end{equation}
Therefore, there exist $\left(  u_{n}^{1\ast},v_{n}^{1\ast}\right)  \in
\frac{1}{n}D_{X^{\ast}\times Y^{\ast}},$ $u_{n}^{2\ast}\in\frac{1}%
{n}D_{X^{\ast}}$ and $v_{n}^{2\ast}\in\frac{1}{n}D_{Y^{\ast}}$ such that
$x_{n}^{1\ast}-u_{n}^{1\ast}\in\widehat{D}^{\ast}F\left(  x_{n}^{1},y_{n}%
^{1}\right)  \left(  v_{n}^{1\ast}-y_{n}^{1\ast}\right)  ,$ $x_{n}^{2\ast
}-u_{n}^{2\ast}\in\widehat{N}\left(  \left(  \overline{x}+\operatorname*{cone}%
L\right)  \cap A,x_{n}^{2}\right)  $ and $y_{n}^{2\ast}-v_{n}^{2\ast}%
\in-\widehat{N}\left(  K,\overline{y}-y_{n}^{2}\right)  \subset K^{+}.$ Using
relation (\ref{1}) we obtain that the sequences $\left(  x_{n}^{1\ast}\right)
,$ $\left(  x_{n}^{2\ast}\right)  ,$ $\left(  y_{n}^{1\ast}\right)  $ and
$\left(  y_{n}^{2\ast}\right)  $ are bounded, and since $X$ and $Y$ are
Asplund spaces, there exist $x_{1}^{\ast}\in X^{\ast},$ $x_{2}^{\ast}\in
X^{\ast}$, $y_{1}^{\ast}\in Y^{\ast}$ and $y_{2}^{\ast}\in Y^{\ast}$ such that
$x_{n}^{1\ast}\overset{w^{\ast}}{\rightarrow}x_{1}^{\ast},$ $x_{n}^{2\ast
}\overset{w^{\ast}}{\rightarrow}x_{2}^{\ast},$ $y_{n}^{1\ast}\overset{w^{\ast
}}{\rightarrow}y_{1}^{\ast},$ $y_{n}^{2\ast}\overset{w^{\ast}}{\rightarrow
}y_{2}^{\ast}$. Obviously, $x_{1}^{\ast}+x_{2}^{\ast}=0$ and $y_{1}^{\ast
}+y_{2}^{\ast}=0.$

Now, if $y_{1}^{\ast}=0,$ then $y_{2}^{\ast}=0$, whence $y_{n}^{2\ast}%
-v_{n}^{2\ast}\overset{w^{\ast}}{\rightarrow}0$ and using the (SNC) assumption
we have that $y_{n}^{2\ast}-v_{n}^{2\ast}\rightarrow0,$ whence $y_{n}^{2\ast
}\rightarrow0,$ so $y_{n}^{1\ast}\rightarrow0.$ Taking into account that $F$
is Lipschitz-like around $\left(  \overline{x},\overline{y}\right)  $ and
using \cite[Theorem 1.43]{Mor-2006-v1}, we obtain that $x_{n}^{1\ast}%
-u_{n}^{1\ast}\rightarrow0$ and since $u_{n}^{1\ast}\rightarrow0,$ we have
that $x_{n}^{1\ast}\rightarrow0.$ Using again (\ref{1}) we obtain that
$x_{n}^{2\ast}\rightarrow0,$ which contradicts the fact that $y_{n}^{2\ast
}\rightarrow0$ and $\left\Vert \left(  x_{n}^{2\ast},y_{n}^{2\ast}\right)
\right\Vert =\frac{1}{2}.$ Hence $y_{1}^{\ast}\neq0.$ Moreover, since
$y_{1}^{\ast}+y_{2}^{\ast}=0,$ $y_{n}^{2\ast}-v_{n}^{2\ast}\overset{w^{\ast}%
}{\rightarrow}y_{2}^{\ast}$, $y_{n}^{2\ast}-v_{n}^{2\ast}\subset K^{+}$ and
$K^{+}$ is weakly-star closed, we obtain that $-y_{1}^{\ast}=y_{2}^{\ast}\in
K^{+}.$

Further, using the hypothesis (iii), for every $n$ large enough, we get that
there exist $l_{n}\in\left(  \overline{x}+\operatorname*{cone}L\right)  \cap
B\left(  x_{n}^{2},\frac{1}{n}\right)  ,a_{n}\in A\cap B\left(  x_{n}%
^{2},\frac{1}{n}\right)  $ such that%
\[
x_{n}^{2\ast}\in\widehat{N}\left(  \left(  \overline{x}+\operatorname*{cone}%
L\right)  \cap A,x_{n}^{2}\right)  \subset\widehat{N}\left(  \overline
{x}+\operatorname*{cone}L,l_{n}\right)  +\widehat{N}\left(  A,a_{n}\right)
+\frac{1}{n}D_{X^{\ast}},
\]
whence, there exist $a_{n}^{\ast}\in\widehat{N}\left(  A,a_{n}\right)  ,$
$l_{n}^{\ast}\in\widehat{N}\left(  \overline{x}+\operatorname*{cone}%
L,l_{n}\right)  $ such that $a_{n}^{\ast}+l_{n}^{\ast}-x_{n}^{2\ast
}\rightarrow0.$ Further, we prove that $\left(  a_{n}^{\ast}\right)  $ or
$\left(  l_{n}^{\ast}\right)  $ is bounded. Suppose by contradiction that both
sequences are unbounded. It follows that for every $n$, there is $k_{n}$
sufficiently large such that%
\begin{equation}
n<\min\left\{  \left\Vert a_{k_{n}}^{\ast}\right\Vert ,\left\Vert l_{k_{n}%
}^{\ast}\right\Vert \right\}  . \label{SB}%
\end{equation}
For simplicity we denote the subsequences $\left(  a_{k_{n}}^{\ast}\right)  $,
$\left(  l_{k_{n}}^{\ast}\right)  $ by $\left(  a_{n}^{\ast}\right)  ,$
$\left(  l_{n}^{\ast}\right)  ,$ respectively. Now, since $a_{n}^{\ast}%
\in\widehat{N}\left(  A,a_{n}\right)  ,$ $l_{n}^{\ast}\in\widehat{N}\left(
\overline{x}+\operatorname*{cone}L,l_{n}\right)  $ we obtain that%
\begin{align*}
\frac{1}{n}a_{n}^{\ast}  &  \in\widehat{N}\left(  A,a_{n}\right)  ,\\
\frac{1}{n}l_{n}^{\ast}  &  \in\widehat{N}\left(  \overline{x}%
+\operatorname*{cone}L,l_{n}\right)  =\widehat{N}\left(  \operatorname*{cone}%
L,l_{n}-\overline{x}\right)  .
\end{align*}
Since
\[
\frac{1}{n}\left\Vert a_{n}^{\ast}+l_{n}^{\ast}\right\Vert \leq\frac{1}%
{n}\left\Vert a_{n}^{\ast}+l_{n}^{\ast}-x_{n}^{2\ast}\right\Vert +\frac{1}%
{n}\left\Vert x_{n}^{2\ast}\right\Vert ,
\]
we obtain that $\frac{1}{n}\left(  a_{n}^{\ast}+l_{n}^{\ast}\right)
\rightarrow0,$ so using again the hypothesis of alliedness we obtain that
$\frac{1}{n}a_{n}^{\ast}\rightarrow0$ and $\frac{1}{n}l_{n}^{\ast}%
\rightarrow0,$ which is in contradiction with relation (\ref{SB}).
Consequently, we obtain that $\left(  a_{n}^{\ast}\right)  ,\left(
l_{n}^{\ast}\right)  \subset X^{\ast}$ are bounded, thus there exist $a^{\ast
},l^{\ast}\in X^{\ast}$ such that $a_{n}^{\ast}\overset{w^{\ast}}{\rightarrow
}a^{\ast}$and $l_{n}^{\ast}\overset{w^{\ast}}{\rightarrow}l^{\ast},$ so
$x_{2}^{\ast}=a^{\ast}+l^{\ast}\in N\left(  A,\overline{x}\right)  +N\left(
\operatorname*{cone}L,0\right)  .$ Now, observe from above that $x_{1}^{\ast
}\in D^{\ast}F\left(  \overline{x},\overline{y}\right)  \left(  y_{2}^{\ast
}\right)  ,$ with $y_{2}^{\ast}\in K^{+}\mathbb{\setminus}\left\{  0\right\}
$ and since $x_{1}^{\ast}+x_{2}^{\ast}=0,$ we get that $0\in D^{\ast}F\left(
\overline{x},\overline{y}\right)  \left(  y_{2}^{\ast}\right)  +N\left(
A,\overline{x}\right)  +N\left(  \operatorname*{cone}L,0\right)  $ with
$y_{2}^{\ast}\in K^{+}\mathbb{\setminus}\left\{  0\right\}  ,$ i.e., the
conclusion.\hfill$\square$\bigskip

We end this section by considering the situation where the objective map is a
single-valued mapping. Consider $f:X\rightarrow\mathbb{R}$ a real-valued
function, take $A\subset X$ and $L\subset S_{X}$ nonempty closed sets. In
order to obtain necessary condition for directional Pareto minimum in the
nonsmooth case, we make use of the penalty function method.

\begin{pr}
Let $\overline{x}\in A$ be a local directional minimum for $f$ on $A$ with
respect to $L$. Suppose that $f$ is Lipschitz continuous around $\overline
{x},$ and $\operatorname*{cone}L$ is convex. In addition, suppose that
$N\left(  A,\overline{x}\right)  \cap\left(  -L^{-}\right)  =\left\{
0\right\}  $ and that either $A$ or $\overline{x}+\operatorname*{cone}L$ is
(SNC) at $\overline{x}$. Then one has%
\[
0\in\partial f\left(  \overline{x}\right)  +N\left(  A,\overline{x}\right)
+L^{-}.
\]

\end{pr}

\noindent\textbf{Proof.} According to the definition of directional minima,
$\overline{x}$ is a local solution of the constrained optimization problem
\begin{equation}
\min f\left(  x\right)  ,\quad x\in\Omega\label{pen1}%
\end{equation}
where $\Omega:=A\cap\left(  \overline{x}+\operatorname*{cone}L\right)  .$
Then, following the well-known Clarke penalization, $\overline{x}$ a solution
of the unconstrained optimization problem%
\[
\min f\left(  x\right)  +kd\left(  x,\Omega\right)  ,\quad x\in X,
\]
where $k>0$ is the Lipschitz modulus of $f.$ By the generalized Fermat rule
and the sum rule for limiting subdifferential, one has%
\begin{align*}
0  &  \in\partial\left(  f+kd\left(  \cdot,\Omega\right)  \right)  \left(
\overline{x}\right)  \subset\partial f\left(  \overline{x}\right)  +k\partial
d\left(  \cdot,\Omega\right)  \left(  \overline{x}\right) \\
&  \subset\partial f\left(  \overline{x}\right)  +N\left(  \Omega,\overline
{x}\right)  .
\end{align*}
Observe that $N\left(  \overline{x}+\operatorname*{cone}L,\overline{x}\right)
=N\left(  \operatorname*{cone}L,0\right)  =L^{-}$ and now we can use
\cite[Corollary 3.5]{Mor-2006-v1} since, according to our assumptions, both
normal qualification condition and the required (SNC) property hold. Then this
allow us to write that%
\[
N\left(  \Omega,\overline{x}\right)  \subset N\left(  A,\overline{x}\right)
+N\left(  \overline{x}+\operatorname*{cone}L,\overline{x}\right)  =N\left(
A,\overline{x}\right)  +L^{-},
\]
and the conclusion follows. \hfill$\square$

\bigskip

Now, we make one step forward by considering the vectorial optimization
problem%
\begin{equation}
\min f\left(  x\right)  ,\quad x\in A, \label{pen3}%
\end{equation}
where $f:X\rightarrow Y$ is a vector-valued function and $A\subset X$ is a
closed set. As before, the ordering cone on $Y$ is $K$.

Consider the following vectorial Lipschitz property for $f$: following
\cite{Ye-NA-2012}, one says that $f$ is $K-$Lipschitz around $\overline{x}\in
X$ of rank $\ell_{f}>0$ if there exist a neighborhood $U$ of $\overline{x}$
and an element $e\in K\cap S_{Y}$ such that for every $x^{\prime}%
,x^{\prime\prime}\in U$%
\[
f\left(  x^{\prime\prime}\right)  -f\left(  x^{\prime}\right)  +\ell
_{f}\left\Vert x^{\prime\prime}-x^{\prime}\right\Vert e\in K.
\]
We record the following result.

\begin{thm}
Let $\overline{x}\in A$ be a local directional Pareto minimum for $f$ on $A$
with respect to $L\subset S_{X}.$ Suppose that:

(i) $f$ is $K-$Lipschitz around $\overline{x}$ of rank $\ell_{f}$ and let $e$
be the element in $K\cap S_{Y}$ given by the Lipschitz property of $f$;

(ii) $K$ is (SNC) at $0;$

(iii) $\operatorname*{cone}L$ is convex, $N\left(  A,\overline{x}\right)
\cap\left(  -L^{-}\right)  =\left\{  0\right\}  $ and that either $A$ or
$\overline{x}+\operatorname*{cone}L$ is (SNC) at $\overline{x}.$

Then for every $\ell>\ell_{f}$, there exist $y^{\ast}\in K^{+}\setminus\{0\}$
and $x^{\ast}\in D^{\ast}f(\overline{x})(y^{\ast})$ such that%
\[
-x^{\ast}\in N(A,\overline{x})+L^{-}\text{ and }\left\Vert x^{\ast}\right\Vert
\leq ly^{\ast}\left(  e\right)  .
\]

\end{thm}

\noindent\textbf{Proof. }Again, directional Pareto minimality of $\overline
{x}$ means that $\overline{x}$ is a Pareto minimum for $f$ on $A\cap
\lbrack\overline{x}+\operatorname*{cone}L].$ We use now a vectorial variant of
Clarke penalization (see \cite[Theorem 3.2 (i)]{Ye-NA-2012}) to deduce that,
for every $\ell>\ell_{f},$ $\overline{x}$ is an unconstrained Pareto minimum
for the function $f(\cdot)+\ell d\left(  \cdot,A\cap\lbrack\overline
{x}+\operatorname*{cone}L]\right)  e.$ We can now use the method from
\cite[Theorem 3.11]{ADS} to deduce that for every $l>l_{f},$ there exist
$y^{\ast}\in K^{+}\setminus\{0\}$, $x^{\ast}\in D^{\ast}f(\overline
{x})(y^{\ast})$ such that%
\[
-x^{\ast}\in N(A\cap\lbrack\overline{x}+\operatorname*{cone}L],\overline{x})
\]
and $\left\Vert x^{\ast}\right\Vert \leq ly^{\ast}\left(  e\right)  .$ Using
again \cite[Corollary 3.5]{Mor-2006-v1}, we have
\[
-x^{\ast}\in N(A,\overline{x})+L^{-},
\]
and this is the conclusion.\hfill$\square$

\section{Pareto directional minima for sets}

As made clear in Definition \ref{df_dir_min_mf} and the subsequent comments,
the notion of directional Pareto minimum is motivated by the case of
(generalized) mappings. However, it is possible to define such a notion for
sets as well. In order to point out this aspect of directional minimality, in
this section we define some appropriate notions and we give, only briefly,
some examples and optimality conditions for them.

Consider, as above, a closed nonempty set $L\subset S_{X}$ and take now $K$ as
a proper closed convex cone in $X.$

\begin{df}
\label{df_dir_min_m}Let $M\subset X$ be a nonempty set. One says that
$\overline{x}\in M$ is a local directional Pareto minimum point for $M$ with
respect to $L$ if
\begin{equation}
(M\cap(\overline{x}+\operatorname*{cone}L)-\overline{x})\cap-K\subset K.
\label{dir_min_m}%
\end{equation}
If $\operatorname*{int}K\neq\emptyset,$ one says that $\overline{x}\in M$ is a
weak directional Pareto minimum for $M$ if
\begin{equation}
(M\cap(\overline{x}+\operatorname*{cone}L)-\overline{x})\cap
-\operatorname*{int}K=\emptyset. \label{dir_wmin_m}%
\end{equation}

\end{df}

It is simple to see that relation (\ref{dir_min_m}) is equivalent to%
\[
(M-\overline{x})\cap\operatorname*{cone}L\cap-K\subset K,
\]
while relation (\ref{dir_wmin_m}) actually means%
\[
(M-\overline{x})\cap\operatorname*{cone}L\cap-\operatorname*{int}K=\emptyset.
\]
Therefore, (\ref{dir_min_m}) is relevant only if $\operatorname*{cone}%
L\cap-K\neq\left\{  0\right\}  ,$ while for (\ref{dir_wmin_m}) it is important
to have $\operatorname*{cone}L\cap-\operatorname*{int}K\neq\emptyset.$

Now, we give an example that justify the above notions of Pareto minimum.

\begin{examp}
Let $\gamma$ be a closed curve described by the following two parametric
equations%
\[
\left\{
\begin{array}
[c]{l}%
x\left(  t\right)  =2+2\cos t\left(  1-\sin t\right) \\
y\left(  t\right)  =\sin t\left(  1-\cos t\right)
\end{array}
\right.  ,\quad\text{ }t\in\lbrack0,2\pi],
\]
$\overline{\gamma}=\operatorname*{int}\gamma\cup\operatorname*{bd}\gamma$ and
the half-plane $H:=\left\{  \left(  x,y\right)  \in\mathbb{R}\times
\mathbb{R}\mid y\geq-x\right\}  .$ Take $K=\mathbb{R}_{+}^{2},$ $\overline
{x}:=(0,0)$ and the directions set $L:=\{\left(  \cos t,\sin t\right)  \mid
t\in\left(  \pi,1.25\pi\right)  \}$. Now, consider $M:=H\cup\left(
\overline{\gamma}\cap-H\right)  $ as a closed subset of $X:=\mathbb{R}^{2}.$
Observing that $\left(  M-\overline{x}\right)  \cap\operatorname*{cone}%
L\cap-K=\left\{  \left(  0,0\right)  \right\}  \subset K\ $and $\left(
M-\overline{x}\right)  \cap-K$ has points that are not in $K\mathbb{\setminus
}\left\{  \left(  0,0\right)  \right\}  ,$ for instance those one that are on
$\gamma$ and have negative $x-$coordinate, we get that $\overline{x}$ is a
local directional Pareto minimum point for $M$ with respect to $L$, but not a
local Pareto minimum point for $M$. Similarly, we have $\left(  M-\overline
{x}\right)  \cap\operatorname*{cone}L\cap-\operatorname*{int}K=\emptyset$ and
$\left(  M-\overline{x}\right)  \cap-\operatorname*{int}K\neq\emptyset,$ so
there exists local weak directional Pareto minimum points, that are not local
weak Pareto minimum points$.$
\end{examp}

In the notation of Definition \ref{df_dir_min_m}, the following optimality
conditions hold.

\begin{thm}
\label{thm_cond_P_d_m}Suppose that $\operatorname*{cone}L\cap
-\operatorname*{int}K\neq\emptyset.$

(i) If $\overline{x}\in M$ is a weak directional Pareto minimum for $M$ with
respect to $L$ then%
\[
T_{B}^{L}(M,\overline{x})\cap-\operatorname*{int}K=\emptyset.
\]

(ii) If for $\overline{x}\in M$ one has%
\[
T_{B}^{L}(\operatorname*{cl}(M+K),\overline{x}))\cap-\operatorname*{int}%
K=\emptyset,
\]
then $\overline{x}$ is a weak directional Pareto minimum for $M$ with respect
to $L.$
\end{thm}

\noindent\textbf{Proof.} (i) Suppose that there exists $u\in T_{B}%
^{L}(M,\overline{x}))\cap-\operatorname*{int}K,$ meaning that $u\in
-\operatorname*{int}K$ and there are $(u_{n})\overset{\operatorname*{cone}%
L}{\longrightarrow}u,(t_{n})\overset{(0,\infty)}{\longrightarrow}0$ such that
for all $n,$ $\overline{x}+t_{n}u_{n}\in M.$ Clearly, for $n$ large enough,
\[
t_{n}u_{n}\in\left(  M-\overline{x}\right)  \cap\operatorname*{cone}%
L\cap-\operatorname*{int}K,
\]
which contradicts the minimality assumption.

(ii) Suppose, again by way of contradiction, that there exists $x\in M$ such
that $x-\overline{x}\in\operatorname*{cone}L\cap-\operatorname*{int}K.$
Consider $(t_{n})\overset{(0,\infty)}{\longrightarrow}0.$ Then, for every $n$
large enough,%
\[
\overline{x}+t_{n}(x-\overline{x})=x+(1-t_{n})(\overline{x}-x)\in
(M+\operatorname*{int}K)\cap(\overline{x}+\operatorname*{cone}L),
\]
whence using the fact that $\operatorname*{cl}(M+\operatorname*{int}%
K)=\operatorname*{cl}(M+K)$ (which, in turn, is easy to prove using the
closedness and the convexity of $K$ which ensures $K=\operatorname*{cl}%
\operatorname*{int}K$) one can write:%
\begin{align*}
-\operatorname*{int}K  &  \ni x-\overline{x}\in T_{B}((M+\operatorname*{int}%
K)\cap(\overline{x}+\operatorname*{cone}L),\overline{x})\\
&  =T_{B}(\operatorname*{cl}\left[  (M+\operatorname*{int}K)\cap(\overline
{x}+\operatorname*{cone}L)\right]  ,\overline{x})\\
&  \subset T_{B}(\operatorname*{cl}\left(  M+\operatorname*{int}K\right)
\cap(\overline{x}+\operatorname*{cone}L),\overline{x})\\
&  =T_{B}(\operatorname*{cl}\left(  M+K\right)  \cap(\overline{x}%
+\operatorname*{cone}L),\overline{x})=T_{B}^{L}(\operatorname*{cl}%
(M+K),\overline{x})),
\end{align*}
and this is in contradiction with the hypothesis.\hfill$\square$

\begin{thm}
Suppose that $\operatorname*{cone}L\cap-K\neq\left\{  0\right\}  .$ If for
$\overline{x}\in M$ one has%
\[
T_{B}^{L}(\operatorname*{cl}(M+K),\overline{x}))\cap-K\subset K,
\]
then $\overline{x}$ is a directional Pareto minimum for $M$ with respect to
$L.$
\end{thm}

\noindent\textbf{Proof.} The proof is similar to that of Theorem
\ref{thm_cond_P_d_m} (ii).\hfill$\square$

\section{Conclusions}

The directional efficiencies introduced in this paper generalize in a
meaningful way the classical situation of Pareto optimality and require
non-trivial adaptations of the usual techniques of investigation used in the
latter case. Besides the results of this paper, we think that our approach
opens new possibilities to model directional situations, especially arising in
vector optimization problems dealing with location issues. We consider that
our concept here introduced is able to capture the situation where some
directions are more important than the others (hence which can be dropped) in
the possible models under consideration. Another possible continuation for
theoretical investigation of directional efficiency is to devise an adapted
(directional) normal limiting cone with respect to a set of directions and to
use it in order to write down more specific optimality conditions for our
concept. All these ideas will be topics for future research.

\bigskip

\noindent\textbf{Acknowledgments:} This research was supported by a grant of
Romanian Ministry of Research and Innovation, CNCS-UEFISCDI, project number
PN-III-P4-ID-PCE-2016-0188, within PNCDI III.

\end{document}